\documentclass[11pt,a4wide,oneside,reqno]{amsart}
\usepackage{amssymb,amsmath,url,a4wide}
\usepackage[utf8]{inputenc}
\usepackage[left=3cm,right=3cm,top=3cm,bottom=2cm]{geometry}

\usepackage[english]{babel}
\usepackage[T1]{fontenc}

\usepackage{stmaryrd}
\usepackage{graphicx}
\usepackage{caption}
\usepackage{lmodern}
\usepackage{textcomp}
\usepackage{mathtools, bm}
\usepackage{amssymb, bm}
\usepackage[unq]{unique}
\usepackage{enumerate}
\usepackage{enumitem}
\newcounter{temp_axiom_enum}
\newcounter{case_enum}
\usepackage{comment}
\usepackage[breaklinks]{hyperref}
\usepackage{parskip}
\usepackage[url=false,doi=true,maxnames=4,giveninits=true,isbn=false]{biblatex}
\usepackage{xpatch}
\xpatchbibmacro{volume+number+eid}{%
  \setunit*{\adddot}%
}{%
}{}{}
\DeclareFieldFormat[article]{number}{\mkbibparens{#1}}
\renewbibmacro{in:}{%
  \ifentrytype{article}{}{\printtext{\bibstring{in}\intitlepunct}}}
\bibliography{bib}

\usepackage{amsthm}

\usepackage{cleveref}
\usepackage{nameref}

\newtheorem{theorem}{Theorem}[section]
\newtheorem{lemma}[theorem]{Lemma}
\newtheorem{observation}[theorem]{Observation}
\newtheorem{conjecture}{Conjecture}[section]
\newtheorem{proposition}[theorem]{Proposition}
\newtheorem{corollary}[theorem]{Corollary}

\newtheorem{question}[theorem]{Question}

\theoremstyle{definition}

\newcommand{\N}{\mathbb{N}}
\newcommand{\lv}{\left\vert}
\newcommand{\rv}{\right\vert}

\title{A proof of the 3/5-conjecture in the domination game}

\author{Leo Versteegen}
\address{Department of Pure Mathematics and Mathematical Statistics, Centre for Mathematical Sciences, Wilberforce Road, Cambridge CB3 0WB, United Kingdom}
\email{lvv23@dpmms.cam.ac.uk}

\date{}

\begin{document}

\maketitle

\begin{abstract}
The \emph{domination game} is an optimization game played by two players, Dominator and Staller, who alternately select vertices in a graph $G$. A vertex is said to be \emph{dominated} if it has been selected or is adjacent to a selected vertex. Each selected vertex must strictly increase the number of dominated vertices at the time of its selection, and the game ends once every vertex in $G$ is dominated. Dominator aims to keep the game as short as possible, while Staller tries to achieve the opposite. In this article, we prove that for any graph $G$ on $n$ vertices, Dominator has a strategy to end the game in at most $3n/5$ moves, which was conjectured by Kinnersley, West and Zamani.
\end{abstract}

\section{Introduction}

Given a graph $G$, the \emph{game domination number} $\gamma_{g}(G)$ is defined as the number of moves the domination game will take if Dominator starts and both players play optimally. The analogous quantity for the version of the game in which Staller starts is denoted by $\gamma_g'$. Since the domination game was introduced by Bre{\v{s}}ar, Klav{\v{z}}ar and Rall in \cite{brevsar_domintaion_initial}, the main objective of its study has been to establish bounds on $\gamma(G)$ and $\gamma_g'(G)$ for different classes of graphs. Typically, these bounds are expressed as a dependency on the number of vertices of $G$.

The most general class of graphs on which the domination is studied, is that of all  graphs without isolated vertices. For this class, Kinnersley, West and Zamani \cite{kinnersley_7/10} conjectured the following.

\begin{conjecture}\label{conj:3/5}
If $G$ is an isolate-free graph of order $n$, then $\gamma_g(G) \leq \frac{3}{5}n$.
\end{conjecture}

There has been some progress towards the verification of \Cref{conj:3/5}. When posing the conjecture, Kinnersley, West and Zamani themselves proved the bound $\gamma_g(G) \leq \lceil 7n/10\rceil$. Bujt{\'a}s \cite{bujtas_2/3_and_min_degree} and Henning and Kinnersley \cite{henning_kinnersley_2/3_and_min_degree} improved upon this independently, by showing $\gamma_g(G) \leq 2n/3$. Subsequently, Bujt{\'a}s improved this further to $\gamma_g(G) \leq 5n/8$ in \cite{bujtas_5/8}. 

For some specific classes of graphs, it has been possible to verify the bound $3n/5$. For instance, it was shown to hold for forests of caterpillars in \cite{kinnersley_7/10}, for forests which do not contain leaves at distance four in \cite{bujtas_forests} and then for a still larger subclass of forests in \cite{schmidt_forests}. Finally, in 2016, Marcus and Peleg \cite{marcus_peleg_forests} proposed a proof for the entire class of forests in a yet unpublished preprint. 

The bound $3n/5$ has also been verified for graphs of minimum degree at least 2 by Henning and Kinnersley in \cite{henning_kinnersley_2/3_and_min_degree} and independently by Bujt{\'a}s in \cite{bujtas_2/3_and_min_degree}. In the same paper \cite{bujtas_2/3_and_min_degree}, Bujt{\'a}s more generally proved upper bounds on $\gamma_g$ for all classes of graphs with a given minimum degree. 

For more bounds for different classes of graphs and information about the the rich ecosystem of variants of the domination game, we refer the reader to the monograph \cite{bookdominationgame}.

In this paper, we prove \Cref{conj:3/5}.

\begin{theorem}\label{thm:main}
If $G$ is an isolate-free graph of order $n$, then $\gamma_g(G)\leq 3n/5$.
\end{theorem}

There are very simple examples showing that \Cref{thm:main} is tight, for instance the path or cycle on five vertices. In our concluding remarks, we will discuss a large family of graphs that are extremal in this sense. It seems possible that no examples beyond this family exist. 

With regard to the Staller-start version of the domination game, we obtain the following corollary from \Cref{thm:main}.

\begin{corollary}\label{cor:Staller-start}
If $G$ is an isolate-free graph of order $n$, then $\gamma_g'(G)\leq 3(n-1)/5 + 1 = 3(n+2)/5$.
\end{corollary}

\subsection*{Acknowledgments}
The author is grateful to be funded by Trinity College of the University of Cambridge through the Trinity External Researcher Studentship. Further, he wishes to express his gratitude to Julien Portier for his valuable comments, to his supervisor Julia Wolf for her generous help in improving this paper, and to Douglas Rall for pointing out an error in the proof of \Cref{lemma:structure}.

\section{Preliminaries}
\subsection{Terminology}
For most of the proof, we will pretend that we are observing a concrete domination game in which we have to instruct Dominator which move to play in a given situation. After one of the players has made a move, we will typically have to update some auxiliary data, which we keep track of throughout the game. Both of these tasks will require some vocabulary, which we introduce here.

Let $G$ be the graph the game is played on. We often abbreviate Dominator and Staller as D and S, respectively. By a \emph{time $t$} in the game we mean a situation in which the combined number of moves played so far by D and S is $t$. Of course, $t$ can never be larger than the total number of moves played at the end of the game. If either player selects a vertex $x$ as their move, we say that they \emph{play $x$}. If a vertex $x$ has a neighbor that has been played or $x$ has been played itself, we say that $x$ is \emph{dominated} and otherwise that $x$ is \emph{undominated}. If a vertex and all its neighbors are dominated, we say that the vertex is \emph{unplayable}.
\subsection{Reductions}
Before we go into the actual proof of \Cref{thm:main}, we will make a few observations that will allow us to simplify the game. The \emph{move-gift} version of the domination game follows the same rules as the standard version, with the following two alterations. Firstly, whenever it is Staller's turn, S may gift a move to D and move directly after D, or gift D an arbitrary number of further moves. Secondly, D may play moves that do not dominate any new vertices. To prove \Cref{thm:main} it is sufficient to show the following result.

\begin{proposition}\label{prop:main}
For each graph $G$ on $n$ vertices, D has a strategy for the (Dominator-start) move-gift version of the domination game on $G$ guaranteeing that the game ends after at most $3n/5$ moves.
\end{proposition}

\begin{proof}[Proof of \Cref{thm:main} from \Cref{prop:main}.]
Of course, S can opt to never gift any moves to D, meaning that Staller's ability to gift moves cannot lead to a worse outcome for S than in the original version. To see that Dominator's ability to play illegal moves can never provide an actual advantage, note that instead of the illegal move, D could simply imagine he played the illegal move, while actually playing a legal move at random. If no legal move is available, the actual game is over earlier than Dominator's imaginary game, allowing us to infer the desired bound.
\end{proof}

Since we can make the same simplifications for the Staller-start version of the game, it is straightforward to infer \Cref{cor:Staller-start} from \Cref{prop:main} as well.

\begin{proof}[Proof of \Cref{cor:Staller-start} from \Cref{prop:main}.]
By removing from $G$ the vertex $x$ that Staller plays as their first move as well as any isolated vertices that arise from removing $x$, we obtain an auxiliary graph $H$. We may now pretend that the players are playing the Dominator-start version of the game on $H$. Because Dominator is allowed to make moves that dominate no new vertices, it does not matter that neighbors of $x$ inside $H$ are already dominated.
\end{proof}

Vertices of degree one play a special role in the domination game and we refer to them as \emph{leaves}, even when the graph under consideration is not a tree. The unique neighbor of a leaf is called the \emph{parent} of the leaf. We point out that if a leaf is dominated, its parent must be dominated as well. We call a parent \emph{dependent} if it has degree 2, and we call the non-leaf neighbor of a dependent parent an \emph{anchor}. We say that the anchor \emph{has} the parent. The following lemma captures an important observation about the structure of graphs that are critical to the conjecture.

\begin{lemma}\label{lemma:structure}
It is enough to prove \Cref{prop:main} for graphs $G$ for which all of the following hold.
\begin{enumerate}[label=(\roman*)]
\item No two leaves are adjacent in $G$, i.e., $G$ contains no isolated edge. 
\item No parent has two leaves.
\item No anchor has three dependent parents.
\end{enumerate}
\end{lemma}
\begin{proof}
Let us assume \Cref{prop:main} is true for graphs for which conditions (i), (ii) and (iii) hold. Let $G$ now be a graph on $n$ vertices that has the minimal number of leaves among all graphs with $\gamma_g(G)>3n/5$. By assumption, there is at least one among the three conditions that $G$ does not meet.

Suppose first that $G$ has an isolated edge $xy$. We consider the graph $H:=G-x-y$. Since $H$ has fewer leaves than $G$, we must have $\gamma_g(H)\leq 3(n-2)/5$. But now D can simply pretend to play the game on $H$ only, pretending further that S is gifting a move when they play one of the vertices $x$ or $y$. If all vertices in $H$ are dominated before $x$ or $y$ have been played and D is to move, D plays $x$ and the game is over. In any case, the game on $G$ lasted at most $3(n-2)/5+1<3n/5$ moves.

Suppose next that $G$ has a parent $x$ with two leaves. In the first move, D can play $x$ and we define $H$ to be the graph obtained from $G$ by removing $x$ and all its leaves. D now pretends that the players are playing the Staller-start game on $H$. By minimality of $G$ and our considerations concerning the Staller-start versions of the game, we know that $\gamma_g'(H)\leq (3(n-3)+2)/5$. Therefore, the game overall takes at most $(3(n-3)+2)/5+1=(3n-2)/5$ moves.

Suppose finally that $G$ has an anchor $u$ with three dependent parents $x,y$ and $z$ that each have a respective leaf $x', y'$ and $z'$. We construct a new graph $H$ by removing $x,y,z,x',y'$ and $z'$ from $G$ and appending a leaf $v$ to $u$. This new graph has at least one leaf less than $G$, and therefore D must have a strategy to end the game on $H$ within $3(n-5)/5$ moves. What is more, we may assume that this strategy does not require D to play $v$, since they may play $u$ instead with no disadvantage.

We instruct D to play according to this strategy and assume first that S also plays only vertices in $H$ until all of them are dominated. At this point, $u$ must have been played to dominate $v$, since neither player can actually play $v$ itself. Therefore, the game is forcedly concluded by playing exactly one vertex from each of the sets $\{x,x'\},\{y,y'\}$, and $\{z,z'\}$.

This leaves the possibility that before the game on $H$ is finished, S plays a vertex from $\{x,y,z,x',y',z'\}$, without loss of generality $x$ or $x'$. In this case, D simply plays $y$ and then continues their auxiliary game on $H$. If S plays $z'$ later on, D treats this as a move-gift in the auxiliary game. Likewise, if S plays $z$ in a situation where it would not be legal to play $v$ in $H$, then D treats this as a move-gift as well, noting that S cannot have previously played played $z'$ or play $z'$ in the future. If S plays $z$ in a situation where $v$ or $u$ are undominated in the auxiliary game, D pretends that S played the vertex $v$. In any case, the game on $G$ lasts at most $3n/5$ moves.
\end{proof}

\section{Core elements of the proof}
\subsection{Colors, points and axioms}
Throughout this section, let $G$ be a graph on $n$ vertices that has no isolated vertices and has the properties (i)-(iii) from \Cref{lemma:structure}. The core idea of the proof, inspired by that in \cite{portiertotaldomination}, is that in most situations, Dominator can play either a vertex that dominates many new vertices or a vertex that makes many vertices unplayable. To see why this is helpful, let $d(t)$ denote the number of dominated vertices at time $t$, and let $u(t)$ denote the number of vertices that are unplayable, either because they have been played already or because they and all their neighbors have otherwise become dominated. If we can prove for some constant $c>0$ and all $t\leq T$ a statement of the form $d(t)+u(t)\geq ct$, then we can infer that for all $t$, we have $t\leq 2n/c$. In particular, the game can last for at most $2n/c$ moves.

To prove \Cref{prop:main}, we would need $c$ to be at least $10/3$. The easiest way to obtain a bound like $d(t)+u(t)\geq 10t/3$ would be by induction, i.e., by showing a statement of the form $d(t+1)-d(t)+u(t+1)-u(t)\geq 10/3$. Unfortunately, this is not possible for various reasons. First of all, there is the problem that even if Dominator has good moves available throughout the game, it is essentially unavoidable that Staller gets to play vertices that increase $d$ and $u$ by only one each. Therefore, the best we can hope for is to find for every $t\geq 0$, a $t'>t$ such that $d(t')-d(t)+u(t')-u(t)\geq 10(t'-t)/3$.

Of greater concern are some longer-term problems. For example, it will typically happen that the game reaches a point where the only vertices that are not yet dominated are leaves. The only way for Dominator to make any progress is to play such a leaf or its parent, thus increasing $d$ by one and $u$ by two.

In \cite{portiertotaldomination}, similar problems occur, and they are dealt with by formulating a fairly large number of linear constraints between the number of moves of certain types. Since some types of moves, many of which occur early on in the game, increase $d+u$ by more than the necessary constant, Dominator can afford to play weaker moves towards the end of the game without falling below the critical threshold overall. In spirit, we take a similar approach here, but to prevent the number of variables and linear constraints from becoming to large, we encode most of the information in terms of colors that we assign to individual vertices. Each color has a fixed number of \emph{points} associated to it as follows.

\begin{itemize}
    \item White: 6 points
    \item Yellow: 5 points
    \item Green: 4 points
    \item Blue: 3 points
    \item Orange: 2 points
    \item Red: 0 points
\end{itemize}

The colors are governed by \emph{coloring axioms}, which describe which colors a vertex is allowed to have depending on how the vertices within distance two are structured and whether they are dominated. Some of these axioms are \emph{permanent} in the sense that they are to be maintained throughout the whole game, while other, \emph{temporary} axioms will be introduced and dropped at the start and end of different phases of the game. Of course, it is integral to the validity of the proof that any axiom holds at the time it is introduced. The permanent axioms are as follows.

\begin{enumerate}[label=\textbf{(P\arabic*)}]
    \item\label{axiom:white} If a vertex is not yet dominated, it is always white. 
    \item\label{axiom:red} Only unplayable vertices may be colored red. 
    \item\label{axiom:green} A vertex adjacent to an undominated leaf or an undominated parent can only be white, yellow or green. 
    \item\label{axiom:yellow} Once a vertex has been colored yellow, it can only be recolored if it is not adjacent to an undominated parent.
    \item\label{axiom:orange} An orange vertex can have at most one undominated neighbor. 
    \item\label{axiom:timing} No vertex is colored yellow before the start of \nameref{phase:parent+2} and no vertex is colored blue or orange before the start of \nameref{phase:big-white-components}. 
    \item\label{axiom:reversal} Once a leaf has been given a color other than white, it may not be colored white again. Furthermore, once a parent has been given a color other than white, it may only be colored white again if its leaf is not white anymore.
\end{enumerate}

Clearly, at the beginning of the game, all vertices must be white. The idea is that whenever we recolor a vertex in a way that decreases the number of points associated with its color, we consider these points to be \emph{earned}. Such earned points are not associated with either player and they are not part of the actual game, but rather a tool we use to keep track of its progress. Although we will sometimes abuse language by saying that a player or a move earns some number of points, it is we, not the players or some automatism, who color vertices.

Tracking the number of points earned up until and including move $t$ by $\pi(t)$, we know that the game is over as soon as $\pi(t)$ reaches $6n$ because all vertices must\footnote{In fact, it is not hard to see that the game must be over once $\pi(t)$ exceeds $6n-10$, but we will not make use of this fact.} be red then. Therefore, to prove \Cref{prop:main}, it is enough to show that D has a strategy such that if the game is not over at time $t$ and $\pi(t)\geq 10t$, then there exists $t'>t$ such that $\pi(t')\geq 10t'$. Mostly, we will take stock of the points earned after every double-move, i.e., after one move by Dominator and one by Staller, so that typically $t'=t+2$, but there may be arbitrarily long sequences that need to be played out before we can validate that enough points have been earned. Occasionally, we will also color a non-white vertex white again, and whenever we do so we will deduct the appropriate number of points from the total number of points earned.

At this point we should remark that the idea of coloring vertices and tracking the number of associated points is not new in the context of the domination game. Rather, it was introduced in \cite{bujtas_forests} and \cite{bujtas_2/3_and_min_degree} and later named the \emph{Bujt{\'a}s discharging method}. Indeed, our approach here is very similar, even to the extent that white, blue and red vertices can largely be thought of in the same way as their counterparts in Bujt{\'a}s' work. An important difference, however, is that in the discharging method, the colors a vertex can have exhaustively and exclusively cover all possible states of a vertex. That is, the graph together with the list of moves up to a given point completely determines the coloring of vertices at that point in the game. In our system on the other hand, there might at any given time be many colorings that satisfy the coloring axioms and would therefore be considered legitimate. 

\subsection{An overview of the different phases}
Our strategy is divided into eleven different phases. To keep the reader from getting lost in minute details, we give a quick overview of these and illustrate why they are necessary. 

To motivate our strategy, it is best to start at the last phase, \nameref{phase:isolates}, which is reached when no two undominated vertices are adjacent and no vertex has two undominated neighbors. If we consider an undominated vertex $x$ during this phase, it must be white by \ref{axiom:white} and the further analysis depends on whether $x$ is a leaf or not. If it is, its parent must be yellow or green by \ref{axiom:green}. If $x$ is not a leaf, $x$ has at least two neighbors neither of which are red by \ref{axiom:red}. In either case, if $x$ or one of its neighbors is played, we can color $x$ and all its neighbors red which earns at least 10 points. Since both players must always play such a move, we can even allow D to choose their moves at random during \nameref{phase:isolates}.

The difficulty lies in playing earlier phases of the game such that \ref{axiom:green} actually holds, that is, ensuring that parents of undominated leaves are at least green. To see why this is difficult, let us rewind to a point where there are still many white vertices left, conceivably arranged in all sorts of configurations. We would like to follow the strategy of Bujt{\'a}s from \cite{bujtas_5/8}, which goes as follows. If no white vertex has three white neighbors, the components of the graph induced by $G$ on the white vertices must all be paths or cycles, and it is then easy enough to exploit this linear structure to earn enough points. If such a vertex $x$ does exist, we would like D to play it, after which we would like to color $x$ red and all its white neighbors blue, which earns 15 points. If in response S plays some vertex $y$, we would like to color $y$ red and any undominated neighbor of $y$ whichever color $y$ previously had, which would earn at least 6 points. Taken together, these two moves would earn a healthy 21 points. Ultimately, we will instruct D to play as described above in \nameref{phase:big-white-components} and \nameref{phase:white-interior}, but one obvious problem is that we might end up coloring a parent blue.\footnote{Another one is that $y$ may be orange, but this turns out to be a minor issue.} Note that universally coloring vertices green instead of blue will not solve this problem. Indeed, if we colored $x$ red and three white neighbors of $x$ green, that would earn an insufficient 13 points. This problem turns out to be so severe that we spend the first seven out of eleven phases making sure that parents and their neighbors are colored in a way that ensures a successful endgame.

The plan is to shield parents from becoming blue by making sure that after \nameref{phase:last-uncles}, every parent is either itself green or has all its non-leaf neighbors colored green or yellow. Unfortunately, this can lead to new problems, especially when dependent parents are involved. Suppose we have a green anchor $u$ whose only neighbors are two dependent parents $x$ and $y$ that are both undominated. It follows that the leaves of $x$ and $y$ are undominated also. If D plays $x$ in such a situation, S can respond by playing $u$, and the best we can do is color $u,x$ as well as the leaf of $x$ red and $y$ green, which earns an insufficient 18 points over two moves. Once such a configuration of vertices is present, there is no way to handle it that earns an average of 10 points per move, and the only solution is to prevent it from arising in the first place. \nameref{phase:anchors} is dedicated to exactly this purpose. There we essentially play one dependent parent from each anchor with two dependent parents (Recall that by \Cref{lemma:structure}, no anchor can have more than two dependent parents).

\begin{figure}[htbp]\centering
    			\includegraphics{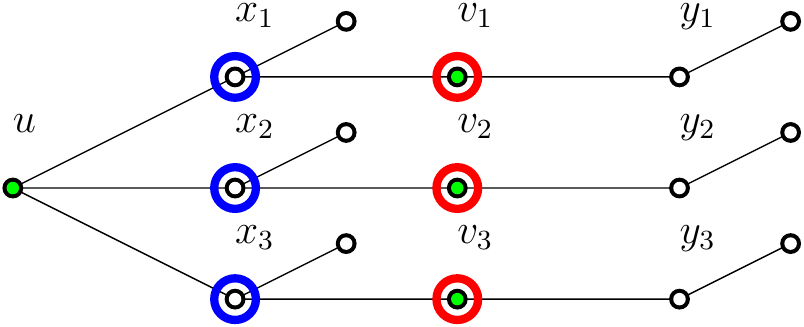}
    			\caption{Moves by Dominator are marked by a blue circle around the played vertex, Staller's moves are red.}
    			\label{fig:trident-danger}
\end{figure}

This still does not solve the problem of anchors in its entirety. Suppose, for example, there was a green vertex $u$ with a neighborhood consisting of three parents $x_1,x_2$ and $x_3$ that are adjacent to distinct anchors $v_1, v_2$ and $v_3$, respectively. Suppose further that each anchor $v_i$ in turn has a dependent parent $y_i$ that is not dominated and that $v_i$ has no other white neighbors. One idea would be to let D play the vertex $x_1$, putting pressure on S to play $v_1$ in order to deny D the very lucrative move $y_1$ at their next turn. Then D could continue by playing first $x_2$ and then $x_3$, to which S might respond by playing $v_2$ and $v_3$. At this point we could color $u,v_1,v_2,v_3,x_1,x_2,x_3$ and all their leaves red as well as $y_1,y_2$ and $y_3$ green (see \Cref{fig:trident-danger}). But sadly, this would only earn 58 points, 2 points short of the required 60 points over six moves. We encourage the reader to verify that D playing $u$ first will earn an insufficient number of points too. Again, the only solution to this problem is to prevent such a configuration from emerging early on, which is the purpose of \nameref{phase:bridges}, \nameref{phase:tridents}, \nameref{phase:uncles+2} and, to some extent, \nameref{phase:parent+2}.

The last major obstacle related to parents arises when two of them are adjacent. For example, suppose that a parent $x$ is not dominated and has two green neighbors $u$ and $v$. Suppose further that $u$ has an undominated leaf and that $v$ is adjacent to a second undominated parent $y$. Although $x$ is not dependent, the undominated leaf of $u$ leads to $x$ behaving similarly to a dependent parent. Namely, if D plays $x$, S can respond by playing $v$. Coloring $x,v$ and the leaf of $x$ red as well as $y$ green would earn 18 points. Note that if the leaf of $u$ did not exist, we could color $u$ red as well for an additional 4 points. During all of the first six phases we need to take great care not to accidentally introduce such a configuration. Furthermore, \nameref{phase:leafed-heads} is dedicated specifically to this cause by giving priority to playing parents that are adjacent to other parents.

Even after we have set up the coloring around parents very carefully in the early phases of the game, it will take  considerable effort to dominate the remaining parents during \nameref{phase:parents}.

\section{The phases in detail}
We are now ready to present the details of Dominator's strategy. We will keep track of the earned points $\pi$, showing that for every time $t$, if the game is not yet finished, there exists $t'>t$ such that $\pi(t')\geq 10t'$. By our earlier remarks, this will prove \Cref{prop:main}. 

The strategy is divided into eleven phases, each of which will described below, including instructions for both moves for D and coloring vertices. The first phase begins with the start of the game and each subsequent phase begins as soon as the previous one has finished. A phase finishes as soon as an \emph{end condition} is met, which we define at the start of each phase. Following the description of each phase, we will typically make some formal observations about the state of the game that will be useful later on. Although these observations are formulated as snapshots of the current situation, it is easy to see that they must remain true for the remainder of them game. 

During the first six phases, which we refer to as \emph{opening phases}, some moves by S will make it necessary for D to play a very specific move in response. We call these responding moves by D \emph{reactive}. If S has gifted a move to D or S has made a move that does not require a specific reaction or it is the first move of the game, D can play what we call an \emph{active} move.

We will include in the end condition of every opening phase that it must be Dominator's turn and that D is not forced to play a reactive move according to their instructions in the current phase. This ensures that D has the first move in every opening phase,\footnote{That is, if the phase takes place at all. Theoretically, the end condition can be met already before D has made the first move in a phase.} and in any such first move, D will earn at least 14 points.

After the first move of a phase, the concept of reactive moves makes it convenient to think of S moving first, and D responding, actively or reactively, in such a way that both moves taken together earn 20 points. If S gifts a move to D, we simply have to make sure that Dominator's move earns at least 10 points. Therefore, if $t$ moves have been played before the start of an opening phase and $\pi(t)\geq 10t$, we can conclude that $\pi(t')\geq 10t'+4$ whenever D has played move number $t'$ during the same phase. We will also show that every move of S which does not require a reaction from D earns at least 6 points by itself. Leading an induction through the opening phases, we can thereby conclude that if $t$ moves have been played before the start of \nameref{phase:last-uncles}, we must have $\pi(t)\geq 10t$. 

From that point onward, we will switch our viewpoint and think of D always moving first. Starting from a time $t$ such that $\pi(t)\geq 10t$, we will always provide D with instructions to advance the game to a time $t'>t$ such that D is to move next and $\pi(t')\geq 10t'$.

During the detailed descriptions of the individual phases, we will typically only give the number of points earned for a particular move without explaining how they combine over several moves to a sufficient total. But now that we have explained the general mechanisms by which the bound $\pi(t)\geq 10t$ is maintained, it will be straightforward for the reader to infer this for themselves.

\subsection*{Phase 1}\label{phase:anchors}
If we have vertices $u,v$ and $x$ such that $v$ is an anchor, $u$ is not a leaf or a dependent parent and $x$ is a parent that is adjacent only to its leaf and the vertices $u$ and $v$, we say that $x$ is an \emph{anchor-bridge} of $u$. Note that if $u$ is also an anchor and $v$ is not a dependent parent, then $x$ is an anchor-bridge of $v$ as well and that no vertex can be both an anchor and an anchor-bridge. A parent is \emph{safe} if it has no undominated neighbor other than possibly its leaf, and a parent is \emph{pseudo-safe} if it is an anchor-bridge, its leaf is white, it has only one other undominated neighbor and this neighbor is an anchor. Finally, we say that a parent is \emph{fully white} if both it and its leaf are white, noting that once a parent ceases to be fully white, it cannot become fully white again by \ref{axiom:reversal}.

This phase ends when it is Dominator's turn, D does not need to play a reactive move and every anchor satisfies one of the following.
\begin{itemize}
\item It is green or red.
\item It has no white dependent parent.
\item It has no white leaf, has one white dependent parent, and it is not adjacent to any other fully white parent.
\end{itemize}

In addition to the permanent axioms, we will maintain the following temporary coloring axioms during the first phase:

\begin{enumerate}[label=\textbf{(T\arabic*)}]
    \item\label{axiom:green-anchors-phase} A vertex can only be green if it is one of the following.
    \begin{itemize}
        \item A safe or pseudo-safe parent.
        \item An anchor that has no white leaf, exactly one fully white dependent parent and at most one undominated anchor-bridge .
    \end{itemize}
    \item\label{axiom:white-dependent-anchor-leaf} If an anchor $u$ has no fully white dependent parent, $u$ has no white leaf and at most one undominated anchor-bridge.
    \setcounter{temp_axiom_enum}{\value{enumi}}
\end{enumerate}

We now give instructions for coloring vertices and moves by D, distinguishing different types of moves by S that might the precede the current move. When we do not mention explicitly a reactive move that D must play, it is understood that D should play one of the active moves given later on.

Suppose first that S has just played a green vertex $x$. If $x$ is a (pseudo-)safe parent with a white leaf, we color $x$ and the leaf red to earn 10 points. If $x$ is an anchor with a white dependent parent, we color $x$ red and its white dependent parent green to earn 6 points. By \ref{axiom:green-anchors-phase}, no other kinds of green vertices exist.

\begin{figure}[htbp]\centering
    			\includegraphics{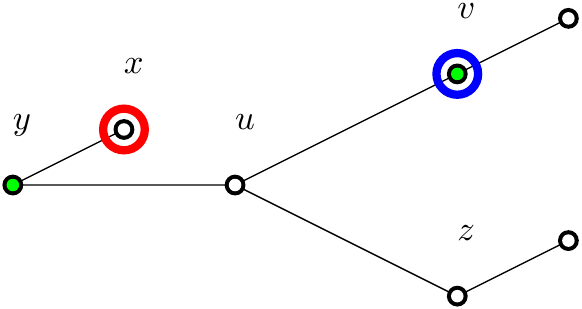}
    			\caption{Dominator's reaction when Staller plays the leaf of a pseudo-safe parent. The picture suggests that $v$ is green, but $v$ might just as well be white. Here and in all subsequent figures, the color configuration shown always leads to the least number of points being earned.}
    			\label{fig:pseudo-safe-defense}
\end{figure}

Next we assume that S has played a white vertex $x$. If $x$ is a leaf of a safe green parent $y$, then we color $x$ and $y$ red to earn at least 10 points. If $y$ is only pseudo-safe, $y$ must be adjacent to a white anchor $u$. If $u$ has a white leaf, D can play $u$ and we color $x,y,u$ and the leaf of $u$ red to earn 22 points. Suppose now that $u$ has no white leaf, and let $z$ be a dependent parent of $u$. If $u$ is adjacent to another parent with a white leaf, we fix such a parent and call it $v$. Then D may play $v$, and we color $v,y$ and their leaves red to earn at least 20 points (see \Cref{fig:pseudo-safe-defense}). If $u$ is not adjacent to any parent with a white leaf other than $z$, we distinguish whether $z$ is fully white or not. If it is, D may play $z$, and we color $y,z$ and their leaves red to earn 22 points. If $z$ is not fully white, the leaf of $z$ is dominated and $z$ must be green or white as $u$ is undominated. Here, D can play $u$ and we color $x,y,z$ and $u$ red to earn at least 20 points.

Suppose next that $x$ is a fully white dependent parent of an anchor $u$. If $u$ has a white leaf or more than one fully white anchor-bridge, D may play $u$ and we color $u,x$, the leaf of $x$ as well as the leaf of $u$ (if it exists) red and all fully white anchor-bridges of $u$ green to earn at least 20 points. If $u$ has no white leaf and at most one fully white anchor-bridge, but $u$ is adjacent to a parent $y$ with a white leaf, $D$ may play $y$, and we color $x,y$ and both their leaves red. We also recolor $u$ white if it was previously green (note that this does not violate \ref{axiom:reversal} since $u$ has no white leaf) and earn a net gain of at least 20 points. Lastly, if $u$ has no leaf and is not adjacent to a second parent with a white leaf, we color $x$ and its leaf red and recolor $u$ white if it was green for a net gain of at least 10 points. If $x$ is the leaf of an undominated dependent parent of an anchor $u$, we proceed in the same way, except for the case where $u$ is undominated, has no white leaf and no second parent with a white leaf, where we only color $x$ red but not its parent. 

If S plays any other white vertex, we simply color that vertex red to earn 6 points. If D does not need to react to the move as described above, D can play an active move. Since the phase is not yet over, there has to exist an anchor $u$ that does not satisfy any of the criteria listed in the end condition. In particular, $u$ must be white and must have a white dependent parent $x$.

If $u$ also has a white leaf, D may play $u$ and we color $u$ and its leaf red and $x$ green to earn 14 points. Likewise, if $u$ has three fully white anchor-bridges or two fully white anchor-bridges and a second white dependent parent, D can still play $u$, and we color it red and all white adjacent anchor-bridges and dependent parents green to earn at least 14 points.

Assuming that none of the above apply to $u$, suppose next that $u$ is adjacent to a fully white parent different from $x$. We now select one such a parent $y$, giving highest priority to dependent parents other than $x$, second highest priority to anchor-bridges and taking any other parent if $u$ has neither of those. We instruct D to play $y$. If the leaf of $x$ is already dominated, $x$ must be dominated as well, and we can color $y$, its leaf and $x$ red to earn 18 points. If the leaf of $x$ is not dominated, we color $y$ and its leaf red and $u$ green to earn 14 points. Note that by our earlier analysis, $u$ has exactly one fully white dependent parent and at most one undominated anchor-bridge after this move.

After the phase is completed, we can make the following observations about the state of the game.

\begin{observation}\label{obs:anchors-no-leaf}
No anchor $u$ has a white leaf.
\end{observation}
\begin{proof}
If $u$ is white, it follows from the end condition and \ref{axiom:white-dependent-anchor-leaf} that $u$ has no white leaf. If $u$ is green, it follows from \ref{axiom:white-dependent-anchor-leaf} that $u$ is not a (pseudo-)safe parent and then from \ref{axiom:green-anchors-phase} that it has no white leaf. 
\end{proof}

\begin{observation}\label{obs:no-two-dependent-parents}
No anchor has two (or more) undominated dependent parents.
\end{observation}
\begin{proof}
If a white anchor had two undominated dependent parents, they would both be fully white and the phase would not be over. If a green anchor had two undominated dependent parents, this would violate \ref{axiom:green-anchors-phase}.
\end{proof}

\begin{observation}\label{obs:anchor-anchor-bridges}
No anchor $u$ has two undominated anchor-bridges.
\end{observation}
\begin{proof}
Suppose $u$ does have two undominated anchor-bridges. By \ref{axiom:red}, $u$ cannot not be red, and by \ref{axiom:green-anchors-phase}, $u$ cannot be green, whence $u$ must be white. Furthermore, by \ref{axiom:white-dependent-anchor-leaf}, $u$ must have a fully white dependent parent, but then $u$ satisfies none of the criteria of the end condition and the phase would not be over.
\end{proof}

\begin{observation}\label{obs:undominated-anchor}
If an anchor $u$ is undominated, $u$ has at most one fully white dependent parent and is not adjacent to any other fully white parent. 
\end{observation}
\begin{proof}
By definition of an anchor, $u$ must have some dependent parent $x$. If $x$ is red, $u$ must be dominated by \ref{axiom:red}. If $x$ is green, $u$ must be dominated by \ref{axiom:green-anchors-phase}. Therefore, $x$ must be white. Since the phase is over, $u$ must satisfy the third criterion of the end condition and therefore cannot be adjacent to another fully white parent.
\end{proof}

All of the observations above must remain true for the rest of the game. To see this, note that the sets of dependent parents, anchors and anchor-bridges are fixed throughout the game, no dominated vertex can become undominated again and no new white leaves or fully white parents can arise due to \ref{axiom:reversal}.

\subsection*{Phase 2}\label{phase:bridges}

A \emph{bridgehead} is a vertex that is not an anchor and has no white leaf, but has a fully white anchor-bridge, and a \emph{trident} is a bridgehead with exactly three fully white anchor-bridges.

This phase ends when it is Dominator's turn, D does not need to play a reactive move and every vertex $u$ satisfies one of the following.
\begin{itemize}
\item $u$ is an anchor.
\item $u$ has no fully white anchor-bridge.
\item $u$ is a dominated bridgehead with at most one (hence, exactly one) fully white anchor-bridge.
\item $u$ is a bridgehead and is adjacent to at most  one (hence, exactly one) fully white parent.
\item $u$ is a dominated bridgehead that is adjacent to exactly two fully white parents.
\item $u$ is a trident.
\end{itemize}

In \nameref{phase:bridges}, we maintain an axiom that is similar to but more permissive than \ref{axiom:green-anchors-phase}.

\begin{enumerate}[label=\textbf{(T\arabic*)}]
    \setcounter{enumi}{\value{temp_axiom_enum}}
    \item\label{axiom:green-bridges-phase} A vertex can only be green if it is one of the following.
    \begin{itemize}
        \item A safe or pseudo-safe parent.
        \item An anchor with a fully white dependent parent.
       \item A bridgehead with at most one fully white anchor-bridge.
       \item A bridgehead that is adjacent to exactly two fully white parents.
    \end{itemize}
    \setcounter{temp_axiom_enum}{\value{enumi}}
\end{enumerate}

We begin by giving instructions on how to react to moves by Staller. If S plays a green bridgehead $u$, we simply color $u$ red and all fully white anchor-bridges of $u$ green to earn at least 6 points. Suppose next that there is a fully white anchor-bridge $x$ with a leaf $y$ such that $x$ is adjacent to a green bridgehead $u$, and S plays $x$ or $y$. By \Cref{obs:undominated-anchor}, the anchor $v\neq u$ that $x$ is adjacent to is dominated. Hence, $x$ is safe, and we may color $x$ and $y$ red and $u$ white for a net gain of 10 points, noting that $u$ cannot have a white leaf by \ref{axiom:green-bridges-phase}. 


In any other case, D reacts and we color as in \nameref{phase:anchors}. If D does not need to play a reactive move, and the phase is not over, then we consider a vertex $u$ that does not satisfy any of the criteria of the end condition.

Suppose first that $u$ has four or more fully white anchor-bridges. By \ref{axiom:red} and \ref{axiom:green-bridges-phase}, $u$ must be white, and D will simply play $u$ after which we color $u$ red and all its adjacent anchor-bridges green to earn 14 points. Likewise, if $u$ has a white leaf and at least one fully white anchor-bridge, $u$ must also be white, and D can play $u$, after which we color it and its leaf red and all fully white anchor-bridges that are adjacent to $u$ green to earn at least 14 points.

Suppose next that $u$ has no white leaf but exactly two fully white anchor-bridges $x$ and $y$. Again, $u$ must be white as otherwise $u$ would either satisfy the fifth criterion or violate \ref{axiom:green-bridges-phase}. Therefore, D can play $y$, and we can color $y$ and its leaf red and $u$ green to earn 14 points. This leaves only the possibility that $u$ is undominated and has only a single fully white anchor-bridge $x$ and some other fully white parent $y$, in which case D can play $y$ just the same.

Before making some observations about the state of the game at the end of this phase, we will have to introduce some more terminology. Let $A\subset V(G)$ be the set of all vertices that are anchors, dependent parents, anchor-bridges or bridgeheads at this point. Among them, we fix the subset $R\subset A$ of vertices that are tridents. Furthermore, let $B$ be the set of vertices that are neither a leaf nor in $A$. The sets $A, B$ and $R$ should be understood as snapshots at the end of \nameref{phase:bridges}, i.e., even if a vertex ceases to be a bridgehead or a trident in the future, we do not remove it from $A$ or $R$, respectively. Note however, that at any time from here on out, all bridgeheads are elements of $A$ by \ref{axiom:reversal}.

\begin{observation}\label{obs:undominated-in-A-adjacent-B}
If an undominated vertex $u\in A$ is adjacent to a parent $x\in B$ with a white leaf, $u$ must be a trident, i.e., $u\in R$.
\end{observation}
\begin{proof}
Clearly, $u$ cannot be a dependent parent or an anchor-bridge as this would place $x$ inside $A$. Furthermore, $x$ must be white by \ref{axiom:red} and \ref{axiom:green-bridges-phase} (recall that  by definition only anchor-bridges can be pseudo-safe), whence by \Cref{obs:undominated-anchor}, $u$ cannot be an anchor. This leaves only the possibility that $u$ is a bridgehead, and by the end condition of \nameref{phase:bridges}, $u$ must be a trident.
\end{proof}

\begin{observation}\label{obs:undominated-parent-in-A}
An undominated parent $x\in A$ is either dependent or an anchor-bridge.
\end{observation}
\begin{proof}
By \Cref{obs:anchors-no-leaf}, $x$ cannot be an anchor and by the end condition of \nameref{phase:bridges}, $x$ is not adjacent to a fully white anchor-bridge.
\end{proof}

\begin{observation}\label{obs:no-cousins-inside-A}
If two adjacent vertices $x$ and $y$ both have a white leaf, and at least one of them is in $A$, both $x$ and $y$ must be dominated.
\end{observation}
\begin{proof}
By \Cref{obs:anchors-no-leaf}, neither $x$ nor $y$ can be an anchor, hence neither of them can be a dependent parent either, and by and the end condition of \nameref{phase:bridges}, they cannot be adjacent to a fully white anchor-bridge. This leaves the possibility that one of them, $x$ say, is itself an anchor-bridge. We then also know that $x$ cannot be white because $y$ is adjacent to $x$. But if $y$ were undominated, $x$ could not be red, by \ref{axiom:red}, or green, by \ref{axiom:green-bridges-phase}, either.
\end{proof}

\subsection*{Phase 3}\label{phase:tridents}
A \emph{residual parent} is a parent that lies in $B$ but has no neighbors in $B$, and if a parent $x\in B$ has only one neighbor in $B$, which we call $z$, we say that $x$ is a \emph{quasi-dependent} parent of the \emph{quasi-anchor} $z$. For a residual or quasi-dependent parent $x$, let $f(x)$ be the number of its neighbors that are tridents. Recall that no vertex has more than three fully white anchor-bridges, and so by \ref{axiom:reversal}, $f$ is decreasing. Finally, a \emph{semi-trident} is a bridgehead with exactly two fully white anchor-bridges.

This phase ends when it is Dominator's turn, D does not need to play a reactive move and no tridents exist anymore. We maintain the following axioms.

\begin{enumerate}[label=\textbf{(T\arabic*)}]
\setcounter{enumi}{\value{temp_axiom_enum}}
    \item\label{axiom:green-trident-phase} A vertex can only be green if it is one of the following.
    \begin{itemize}
        \item A safe or pseudo-safe parent.
        \item An anchor with a fully white dependent parent.
       \item A bridgehead with at most one fully white anchor-bridge.
        \item A semi-trident.
        \item A quasi-dependent parent $x$ with $f(x)=0$ and a white leaf whose quasi-anchor has a fully white quasi-dependent parent $y$ with $f(y)=0$.
    \end{itemize}
    \item\label{axiom:semi-trident-domination} If a vertex $u\in R$ is not a trident anymore, then it must be dominated.
    \item\label{axiom:trident-count-residual} A semi-trident cannot be adjacent to an undominated residual parent $x$ with $f(x)=0$.
    \item\label{axiom:trident-count-transitional} If a semi-trident is adjacent to an undominated quasi-dependent parent $x$ with $f(x)=0$, the quasi-anchor $u$ of $x$ can have no other undominated quasi-dependent parent $y$ with $f(y)=0$.
    \setcounter{temp_axiom_enum}{\value{enumi}}
\end{enumerate}

Note that \ref{axiom:trident-count-residual} and \ref{axiom:trident-count-transitional} are satisfied at the start of the current phase because by the end condition of \nameref{phase:bridges}, no semi-trident is adjacent to any undominated parent other than its two anchor-bridges. Note further 
that a new semi-trident can only arise if the anchor-bridge or leaf of an anchor-bridge of a trident is given a new color.

We first give instructions for when S plays a trident, a fully white anchor-bridge of a trident or its leaf, a quasi-dependent parent $x$ with $f(x)=0$ that is not safe, or the leaf of such a parent. For all other moves by S, we proceed as in \nameref{phase:bridges}, including for moves in which S plays semi-tridents, which are handled just as any other bridgehead. We emphasize that none of the reactions from \nameref{phase:bridges} (including those taken from \nameref{phase:anchors}) involve playing or coloring any fully white anchor-bridge of a trident or the leaf of such an anchor-bridge, which is crucial to avoiding violations of the temporary axioms of \nameref{phase:tridents}.

If S plays a trident $u$, we can simply color $u$ red and its anchor-bridges green, noting that \ref{axiom:trident-count-residual} and \ref{axiom:trident-count-transitional} respectively are unaffected since every parent adjacent to $u$ must now be dominated. Suppose next that S plays a vertex $x$ that is a fully white anchor-bridge of a trident $u$ or a leaf such an anchor-bridge $x$. In this case, D can play $u$ themselves, and we color $u$, $x$, the leaf of $x$ red and the remaining two fully white anchor-bridges of $u$ green to earn 22 points.

If S plays a green quasi-dependent parent $x$ with a white leaf, we simply color $x$ and the leaf red for 10 points. If S plays instead the leaf of a green quasi-dependent parent $x$, and $x$ is not safe, we note that by \ref{axiom:green-trident-phase}, since $x\notin A$, the quasi-anchor $u$ of $x$ must have a fully white quasi-dependent parent $y$ with $f(y)=0$. Now D will play $y$, and we will color $x,y$ and both their leaves red to earn 22 points. Note that by \Cref{obs:undominated-in-A-adjacent-B} and \ref{axiom:semi-trident-domination}, all other potential green quasi-dependent parents $z$ with $f(z)=0$ and adjacent to $u$ are now safe, maintaining \ref{axiom:green-trident-phase}. For the same reason, if S plays the white vertex $y$ themselves, we can color $y$ and its leaf red for 12 points without violating the axiom.

If D is not required to play a reactive move, D should, whenever possible, play one of the three fully white anchor-bridge $z_1$ of a trident $u$, and we subsequently color $z_1$ as well as the leaf of $z_1$ red and $u$ green. By \ref{axiom:green-trident-phase}, $u$ must have been white before so that the suggested coloring would earn 14 points, but since $u$ would become a semi-trident, it might violate \ref{axiom:trident-count-residual} or \ref{axiom:trident-count-transitional}. We now analyze how these violations might occur and describe what other moves to play instead.

Suppose first that playing $z_1$ would lead to a violation of \ref{axiom:trident-count-residual}. Then there is a trident $u$ adjacent to a residual parent $x$ with $f(x)=1$. In this case, D can play $u$, and we color $u$ red and its anchor-bridges as well as $x$ green. Note that by \Cref{obs:undominated-in-A-adjacent-B} and \ref{axiom:semi-trident-domination}, $x$ can have no undominated neighbor other than its leaf. 

\begin{figure}[htbp]\centering
    			\includegraphics{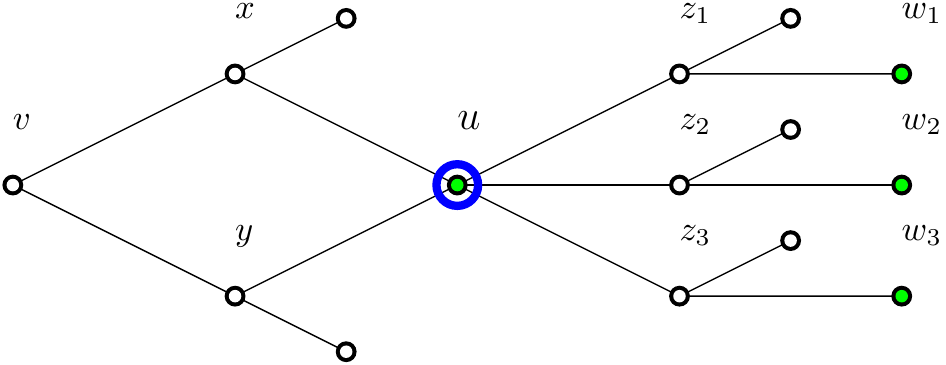}
    			\caption{The trident $u$ is adjacent to the two quasi-dependent parents $x$ and $y$. The anchor-bridges $z_1,z_2$ and $z_3$ are each adjacent to their own anchors $w_1, w_2$ and $w_3$.}
    			\label{fig:trident}
\end{figure}

If playing $z_1$ would instead lead to a violation of \ref{axiom:trident-count-transitional}, then there must be two undominated quasi-dependent parents $x$ and $y$ with a shared quasi-anchor $v$ such that either $f(x)=1$, $f(y)=0$ and $x$ is adjacent to $u$ or $f(x)=f(y)=1$ and both $x$ and $y$ are adjacent to $u$ (see \Cref{fig:trident}). In the former case, D can earn 14 points by playing $u$, coloring $u$ red and its anchor-bridges as well as $x$ green. Note that $x$ being green does not violate \ref{axiom:green-trident-phase} due to the existence of $y$. In the latter case, D can still play $u$, color $u$ red and its anchor-bridges as well as $x$ (but not $y$!) green.

Let us now make some observations about the state of the game at the end of \nameref{phase:tridents}. To see that they will remain true for the rest of the game, observe that no new semi-tridents can arise since no vertex is adjacent to more than two fully white anchor-bridges.

\begin{observation}\label{obs:semi-trident-transitional} If a quasi-anchor $u$ has an undominated quasi-dependent parent $x$ such that $x$ is adjacent to a semi-trident, then $u$ can have no second undominated quasi-dependent parent. 
\end{observation}
\begin{proof}
This is a direct consequence of \ref{axiom:trident-count-transitional} and the fact that $f\equiv 0$.
\end{proof}

\begin{observation}\label{obs:semi-trident-residual}
A semi-trident cannot be adjacent to an undominated residual parent.
\end{observation}
\begin{proof}
This follows from \ref{axiom:trident-count-residual} and the fact that $f\equiv 0$.
\end{proof}

\begin{observation}\label{obs:undominated-in-A-adjacent-B-full}
No undominated vertex $u\in A$ is adjacent to a parent $x\in B$ with a white leaf.
\end{observation}
\begin{proof}
Follows from \Cref{obs:undominated-in-A-adjacent-B} and \ref{axiom:semi-trident-domination}.
\end{proof}

\begin{observation}\label{obs:anchor-bridge-uncles}
Any vertex $u$ that has a fully white anchor-bridge has no white leaf and satisfies at least one of the following.
\begin{itemize}
\item $u$ is dominated and has at most (hence, exactly) one undominated anchor-bridge.
\item $u$ is a dominated semi-trident.
\item $u$ is adjacent to at most (hence, exactly) one fully white parent and $u\notin R$.
\end{itemize}
\end{observation}
\begin{proof}
If $u$ has a white leaf or does not satisfy any of the criteria above, then it must have been a trident at the at the end of \nameref{phase:bridges} (to see that $u$ cannot have been an anchor, recall \Cref{obs:anchors-no-leaf}, \Cref{obs:anchor-anchor-bridges} and \Cref{obs:undominated-anchor}). Hence $u$ must be in $R$. In particular, $u$ is dominated and has no white leaf, and since no more tridents exist, $u$ must have become a semi-trident or a vertex with at most one undominated anchor-bridge. 
\end{proof}

\subsection*{Phase 4}\label{phase:leafed-heads}
This phase ends when it is Dominator's turn, D does not need to play a reactive move and no white quasi-anchor has both a white leaf and a white quasi-dependent parent. We maintain the following temporary axiom, which is essentially \ref{axiom:green-trident-phase} from \nameref{phase:tridents} reformulated to accommodate the fact that $f\equiv 0$.

\begin{enumerate}[label=\textbf{(T\arabic*)}]
\setcounter{enumi}{\value{temp_axiom_enum}}
    \item\label{axiom:green-parent+2-phase} A vertex can only be green if it is one of the following.
    \begin{itemize}
        \item A safe or pseudo-safe parent.
        \item An anchor with a fully white dependent parent.
        \item A bridgehead.
        \item A quasi-dependent parent with a white leaf whose quasi-anchor has a fully white quasi-dependent parent.
    \end{itemize}
    \setcounter{temp_axiom_enum}{\value{enumi}}
\end{enumerate}

D reacts where necessary to Staller's move as during \nameref{phase:tridents}. If no reaction is needed while the phase is ongoing, D plays a white quasi-anchor $u$ that has both a white quasi-dependent parent and a white leaf. We color $u$ and its leaf red and all quasi-dependent neighbors of $u$ green, noting that they are safe by \Cref{obs:undominated-in-A-adjacent-B-full}. This earns 14 points.

\begin{observation}\label{obs:quasi-anchor-no-leaf}
Let $x$ be a quasi-anchor with a white leaf and let $y$ be a quasi-dependent parent of $x$ with a white leaf. Then both $x$ and $y$ must be dominated.
\end{observation}
\begin{proof}
If $x$ is undominated, $x$ cannot have any fully white quasi-dependent parent as otherwise \nameref{phase:leafed-heads} would not be over. In particular, $y$ cannot be white, and by \ref{axiom:green-parent+2-phase}, $y$ cannot green either. Hence, $x$ must be dominated.

If $y$ is undominated, then in turn $x$ cannot be white since the phase is over. Furthermore, \ref{axiom:green-parent+2-phase} precludes $x$ from being green since $x$ is not in $A$, $x$ is not safe and $x$ is not a green quasi-dependent parent of $y$ either, as $y$ has no other neighbor in $B$ than $x$.
\end{proof}

\subsection*{Phase 5}\label{phase:parent+2}
This phase ends when it is Dominator's turn, D does not need to play a reactive move and every fully white parent has at most two white neighbors (one excluding its leaf).  We will continue to maintain \ref{axiom:green-parent+2-phase} during this phase.

From this phase onwards, we may color vertices yellow without violating \ref{axiom:timing}. Whenever S plays a yellow vertex, we color this vertex red and all its neighbors yellow to earn at least 6 points. As for other moves by S, we color vertices and instruct D to react where necessary as in \nameref{phase:tridents} and \nameref{phase:leafed-heads}.  When no reactive move is needed, D simply plays a fully white parent with at least two other white neighbors. We then color parent and leaf red and the other neighbors yellow, earning 14 points.

\begin{observation}\label{obs:core-parents}
An fully white parent $x$ can have at most one neighbor in $B$ that is not yellow. 
\end{observation}
\begin{proof}
By \ref{axiom:red}, $x$ cannot have a red neighbor and by the end condition of the current phase, $x$ can have at most one white neighbor besides its leaf. Suppose $x$ has a green neighbor $y\in B$. By  \ref{axiom:green-parent+2-phase}, $y$ would have to be a quasi-dependent parent with a white leaf, but $x$ cannot be a quasi-anchor by \Cref{obs:quasi-anchor-no-leaf} and $x$ cannot be in $A$ by \Cref{obs:undominated-in-A-adjacent-B-full}.
\end{proof}

We remark that this observation will also always remain true, due to \ref{axiom:yellow}. 

\subsection*{Phase 6}\label{phase:uncles+2}
This phase ends when it is Dominator's turn, D does not need to play a reactive move and no white vertex is adjacent to two or more fully white parents or has a white leaf and is adjacent to one fully white parent (no white vertex can have a white leaf and two white neighbors after \nameref{phase:parent+2}). We say that a quasi-anchor is \emph{mixed} if it is adjacent to a quasi-dependent parent and a non-quasi-dependent parent that are both  undominated.

We maintain the following axiom, which is a further relaxation of \ref{axiom:green-parent+2-phase}.

\begin{enumerate}[label=\textbf{(T\arabic*)}]
\setcounter{enumi}{\value{temp_axiom_enum}}
    \item\label{axiom:green-uncles+2-phase} A vertex can only be green if it is one of the following.
    \begin{itemize}
       \item A safe or pseudo-safe parent.
        \item An anchor with a fully white dependent parent.
        \item A bridgehead.
        \item A quasi-dependent parent with a white leaf whose quasi-anchor has a fully white quasi-dependent parent.    
        \item A vertex in $B$, which
        \begin{itemize}
        	\item has no white leaf,
        	\item is adjacent to a white safe parent that has no other green neighbors,
        	\item is not adjacent to more than two fully white parents overall and
        	\item is not a mixed quasi-anchor.
        \end{itemize}
    \end{itemize}
    \setcounter{temp_axiom_enum}{\value{enumi}}
\end{enumerate}

In this phase, whenever S plays a yellow vertex $u$, we fix an undominated neighbor $x$ of $u$. If $x$ is a safe parent and has a green neighbor $v\in B$, we color $u$ red, $v$ yellow and $x$ green to make a net gain of 6 points. If $x$ is not a safe parent or does not have a green neighbor in $B$, we simply color $u$ red and $x$ yellow to earn 6 points. Furthermore, if $S$ plays a white safe parent $x\in B$ or the leaf of such a parent, we note that $x$ must be undominated and color both $x$ and its leaf red and, if it exists, the unique green neighbor of $x$ in $B$ yellow for a net gain of 11 points. For any other move that Staller may play, we react as in \nameref{phase:uncles+2}.

If no reactive move is needed, we first color all vertices that are white and dominated yellow before we search for an active move. Suppose first that there exists a white (hence undominated) parent $x$ adjacent to another white parent $y$. By the end condition of \nameref{phase:parent+2}, $y$ can have no undominated neighbor other than its leaf. Also, $y$ cannot yet be a safe parent since $x$ is white. Therefore, D can play $x$, and we may color $x$ and its leaf red and $y$ green to earn 14 points, without violating \ref{axiom:green-uncles+2-phase}.

Suppose next that there exists a white vertex $x$ that is adjacent to four or more white parents. These parents must be fully white and are not safe but have no undominated neighbors other than $x$ or their leaf. Hence, D can play $x$, and we can color $x$ red and all white parents adjacent to it green to earn at least 14 points. 

This leaves the possibility that there exists a white vertex $u\in B$ without a leaf but with two or three adjacent fully white parents. We could let D play one of them and then color it as well as its leaf red and $u$ green to earn 14 points, but we will want to give special attention to which one of the parents to play if $u$ is a mixed quasi-anchor. Note that either at most one of the fully white parents adjacent to $u$ is quasi-dependent or at most one of them is not quasi-dependent. By playing the kind of which there is only one and then coloring as suggested above, D can achieve either that $u$ is not adjacent to a fully white quasi-dependent parent or that every fully white parent adjacent to $u$ is quasi-dependent. Of course, if $u$ is not a mixed quasi-anchor to begin with, D can play any one of its fully white parents safely.

\begin{observation}\label{obs:yellow-bad-uncle} If a vertex in $B$ is a mixed quasi-anchor or is adjacent to three fully white parents, then it must be yellow.
\end{observation}
\begin{proof}
It is clear that such a vertex cannot be red or white, and that it is not green follows from \ref{axiom:green-uncles+2-phase}. Indeed, if a vertex is a safe parent, then it cannot have any undominated adjacent parent. A vertex in $B$ cannot be an anchor or have an undominated anchor-bridge, and a quasi-dependent parent with a white leaf cannot be adjacent to any undominated parents by \Cref{obs:undominated-in-A-adjacent-B-full} and \Cref{obs:quasi-anchor-no-leaf}.
\end{proof}

\begin{observation}\label{obs:yellow-cousin}
If two parents with white leaves are neighbors, at least one of them is dominated. Furthermore, if one of them is not dominated, then the other one is yellow.
\end{observation}
\begin{proof}
The first statement follows immediately from the end condition of \nameref{phase:uncles+2}. For the second statement, let $x$ and $y$ be two adjacent parents such that $y$ and the leaves of $x$ and $y$ are undominated. We have to show that $x$ is not green. Suppose towards a contradiction that it is. We show that none of the criteria of \ref{axiom:green-uncles+2-phase} apply to $x$. By \Cref{obs:no-cousins-inside-A}, neither $x$ nor $y$ can be in $A$, $x$ is not safe and by \Cref{obs:quasi-anchor-no-leaf}, $x$ cannot be a quasi-dependent parent. Of course, $x$ has a white leaf, so the last criterion does not apply either.
\end{proof}

\subsection*{Phase 7}\label{phase:last-uncles}
This phase lasts until for each parent, either the parent itself or all of its non-leaf neighbors are dominated. 

We keep no temporary axiom in this phase, and D does not need to play any reactive moves. Whenever S plays a vertex $x$, we simply color the vertex red and all white neighbors of $x$ in the color that $x$ had before, which earns at least 6 points. If the phase is not yet finished, then there must be a white parent with a white leaf and another white neighbor. When it is Dominator's turn, he plays such a parent. We then color parent and leaf red and the other white vertex that is adjacent to the parent green, earning 14 points. 

\subsection*{Phase 8}\label{phase:big-white-components}
This phase lasts until the \emph{undominated graph}, i.e., the subgraph of $G$ induced by the set of undominated vertices, consists only of isolated vertices and edges. 

From this phase onwards, vertices may be colored orange and blue. We remark that if either player plays an undominated vertex $v$, we never violate \ref{axiom:green} by coloring the undominated neighbors of $v$ blue or orange as by the end condition of \nameref{phase:last-uncles}, an undominated neighbor of $v$ cannot be adjacent to an undominated leaf or parent.  During \nameref{phase:big-white-components}, we maintain the following additional axiom.

\begin{enumerate}[label=\textbf{(T\arabic*)}]
\setcounter{enumi}{\value{temp_axiom_enum}}
\item\label{axiom:no-orange-white-white} If an orange vertex has an undominated neighbor, this neighbor may not have an undominated neighbor itself.
\setcounter{temp_axiom_enum}{\value{enumi}}
\end{enumerate}

As before, we can guarantee that every move by S earns at least 6 points: whichever vertex $x$ Staller plays, we color $x$ red and all its undominated neighbors in the color that $x$ had previously. Note that neither \ref{axiom:green} nor \ref{axiom:orange} are violated due to the end condition of \nameref{phase:last-uncles} and \ref{axiom:no-orange-white-white}.

If there exists an undominated vertex $x$ with three or more undominated neighbors, D may play $x$. By coloring $x$ red and all its undominated neighbors blue, at least 15 points are earned, which gives at least 21 points when taken together with the move of S.

If no undominated vertex has three undominated neighbors, we may assume that all remaining components of the undominated graph are paths or cycles. If one of these components is a path on at least three vertices, let us consider the three vertices $x,y,z$ at one of its ends so that $x$ has no undominated neighbors except $y$. After D plays $y$, we color $x$ and $y$ red and $z$ blue, which earns 15 points.

This leaves the possibility that there exists a component that is a cycle $x_1,\ldots,x_k$. If $k\geq 7$, D plays $x_1$. Suppose first that S responds by playing the vertex $x_k$. Then D continues with the move $x_4$, upon which we color the vertices $x_1,x_2,x_3,x_4$ and $x_k$ red and as $x_{k-1}$ as well as $x_5$ blue. Together with the 6 points from the move S plays before it is Dominator's turn again, this gives 42 points. It is easy to see that the same analysis applies not only if S replies by playing $x_2$ rather than $x_k$, but also if they play a vertex outside of the component under consideration. Should S play instead a vertex $x_i$ with $i\notin\{2,k\}$, things become even simpler. We may assume without loss of generality that $x_{i+1}$ was also white, so that by coloring $x_1$ and $x_i$ red as well as $x_k,x_2$ and $x_{i+1}$ blue, overall 21 points are earned. 

For $k=3$ and $k=6$, it is fairly easy to work out a strategy for D, so we leave these cases to the reader and consider only $k=4$ and $k=5$. For $k=4$, D simply plays $x_1$, colors $x_1$ red and both $x_2$ and $x_4$ orange, earning 14 points. For $k=5$, after D plays $x_1$, S can either play outside of the component, one of $x_3$ and $x_4$, or one of $x_2$ and $x_5$. In the first case, D continues by playing $x_3$ and we color all vertices of the cycle red. Together with the two moves by Staller, 42 points are earned this way. In the second case, we can color the whole cycle red immediately to earn 30 points from two moves. In the third case, we may assume by symmetry that S played $x_2$, and we color $x_1$ as well as $x_2$ red and $x_3$ as well as $x_5$ orange, earning 20 points.

\subsection*{Phase 9}\label{phase:white-interior}
This phase lasts until every white vertex is either a leaf, a parent or has only orange neighbors.

We no longer adhere to \ref{axiom:no-orange-white-white}. That is, if an orange vertex $x$ has an undominated neighbor $y$, $y$ may have another undominated neighbor $z$. If in such a configuration, Staller plays $x$, we will color $x$ red and $y$ orange, earning 6 points without infringing on the remaining axioms.

Below is a list of several cases in which D has a move that earns at least 14 points. Having found the first applicable case on the list, D plays the move suggested for this case and we color the vertices appropriately. Once this is done, the idea is to loop back to start this process anew until none of the cases applies any longer. However, at the start of each run of the loop, we will recolor vertices in the color associated with the least number of points that does not violate the coloring axioms. For example, when a dominated vertex that has been colored blue does not actually have any undominated neighbor anymore, we will color it red. This allows us to avoid redundant case distinctions, as it has the following consequences.

\begin{itemize}
    \item No dominated vertex is white.
    \item Every unplayable vertex is red.
    \item Every green vertex is adjacent to either a white parent or a white leaf.
    \item Every blue vertex has at least two white neighbors.
    \item Every yellow vertex is adjacent to a white parent.
\end{itemize}

We remark that by the first point, the properties of being undominated, white or fully white are now synonymous for any parent. In particular, by the end condition of \nameref{phase:big-white-components} no white vertex has two white neighbors.

For the purpose of brevity, we will say that a vertex is \emph{citrus} if it is yellow or green. Furthermore, we will typically only state which vertex D should play and how many points can be earned by this move. We will then leave it to the reader to work out how to color the involved vertices to earn (at least) that number of points.

\begin{enumerate}[label=Case \arabic*., ref=Case \arabic*]
    \item\label{case:endangered-vertex} If there exist vertices $u,x$ and $y$ such that $x$ and $y$ are white neighbors and the white neighborhood of $u$ is a non-empty subset of $\{x,y\}$, then D can earn at least 14 points by playing $x$.

    \item \label{case:triple-white} If there exists a blue vertex $u$ with three white neighbors, then D can play $u$, and we color $u$ red and all white neighbors of $u$ orange, earning 15 points.
    
    \item\label{case:bwwb} If there exists an induced path $uxyv$ such that $x$ and $y$ are white and $u$ and $v$ are blue, then both $u$ and $v$ have exactly two white neighbors by exclusion of \ref{case:triple-white}, and D can earn 14 points by playing $x$.
    
    \item\label{case:wbwb} If there exists a path $uxv$ such that $x$ is white, $u$ is blue and $v$ is blue or orange, then D can earn at least 15 points by playing $u$ if $x$ does not have a white neighbor (note that $u$ must have a second white neighbor) or 14 points by playing $x$ if $x$ does have a white neighbor.
    
    \item\label{case:parent} If there exists a citrus parent $x$ with a white leaf and another white neighbor that is not itself a parent, then D can earn 14 points by playing $x$.
    
    \item\label{case:single-tail} If there exists a citrus vertex that has no white leaf, is adjacent to exactly one white parent and has at most one additional white neighbor, then D can earn 14 points by playing the white parent.

    \item\label{case:grandparent+2} If there exists a citrus vertex $x$ with at least two white non-parent neighbors, then D can earn 14 points by playing $x$ (note that $x$ must have a white leaf or be adjacent to a white parent).

    \item\label{case:triple-tail} If there exists a citrus vertex $x$ that is adjacent to at least three white parents and has an additional white neighbor that is not a parent, then D can earn 14 points by playing $x$.
    \item\label{case:double-tail+} If there exists a citrus vertex $x$ with at least three white neighbors at least one of which does not have a white neighbor itself, then D can earn 14 points by playing $x$.
    \setcounter{case_enum}{\value{enumi}}
\end{enumerate}

We include in the loop one additional case that may require a sequence of four moves (two per player) being played out before a sufficient average number of points has been earned.

\begin{enumerate}[label=Case \arabic*., ref=Case \arabic*]
\setcounter{enumi}{\value{case_enum}}
\item\label{case:zugzwang}
    Suppose there exists a white edge $xy$ such that for every vertex $z$ one of the following is true.
    \begin{itemize}
        \item $z$ has two white neighbors outside of $\{x,y\}$,
        \item $z$ has a white leaf, or
        \item $z$ is red or white.
    \end{itemize}
D starts by playing $x$, and suppose that S replies by playing some vertex $z$. We assume first that $z$ is not a (white) leaf, and color $x$ and $y$ red to earn 12 points and show that we can earn at least 8 more points from Staller's move. 

If $z$ has two white neighbors outside of $\{x,y\}$, if $z$ has a white leaf, or if $z$ is itself white with a white neighbor (neither $x$ nor $y$ can be this white neighbor!), it is easy to see that this is possible. Also, if $z$ is white but does not have a white neighbor and is not a leaf, it cannot not have a citrus neighbor either by exclusion of \ref{case:parent}, \ref{case:single-tail} and \ref{case:double-tail+}. If $z$ is white and has two blue neighbors, we can color $z$ red and its blue neighbors orange (by exclusion of \ref{case:triple-white}), earning 8 points overall. This leaves the possibility that $z$ has an orange neighbor, in which case we can color it red alongside $z$ to earn 8 points.
    
Suppose now that $z$ is a leaf with a citrus parent $u$. We may assume that $u$ is adjacent to a white parent $v$, since otherwise, by exclusion of \ref{case:parent}, we could color $x,y,z$ and $u$ red to earn 22 points. By exclusion of \ref{case:double-tail+}, $u$ has no white neighbors other than $z$ and $v$. But then D can play $v$ in the next move, and we can color $x,y,z,u,v$ and the leaf of $v$ red to earn at least 34 points. Adding the minimum number of 6 points from Staller's second move, we ascertain that at least 40 points have been earned over four moves. 
\setcounter{case_enum}{\value{enumi}}
\end{enumerate}

If none of the above cases apply and there is still a white vertex $x$ which is not a leaf or a parent but has a non-orange neighbor, then we can deduce it must lie in a very specific configuration. Once we have done so, we describe which vertices D should play in this configuration and then go back to the start of the loop.

Note first that $x$ must have a white neighbor. Indeed, if $x$ did not have a white neighbor, then it could not have a citrus neighbor either by exclusion of \ref{case:parent}, \ref{case:single-tail} and \ref{case:double-tail+}. Neither can $x$ have only blue and orange neighbors by exclusion of \ref{case:wbwb}. We denote the white neighbor of $x$ by $y$ and observe that all additional neighbors of $x$ and $y$ must be blue or citrus by exclusion of \ref{case:endangered-vertex}. By exclusion of \ref{case:endangered-vertex} and \ref{case:bwwb}, they cannot both have a blue neighbor, but neither of them is a leaf, so we may assume without loss of generality that $y$ has a citrus neighbor $u$.  

\begin{figure}[htbp]\centering
    			\includegraphics{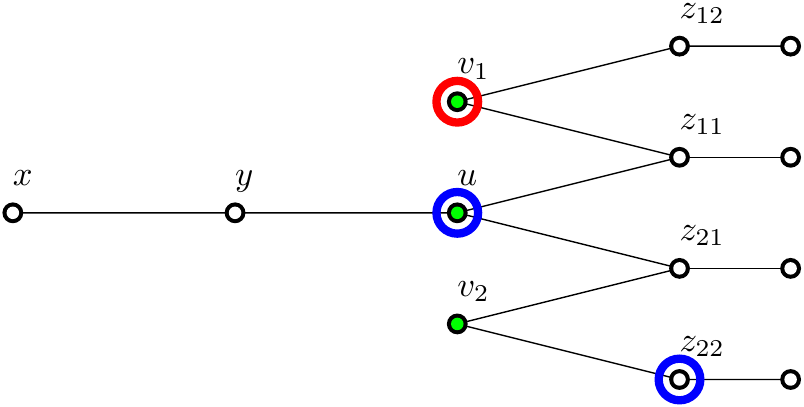}
    			\caption{The situation in which $z_{11},z_{12},z_{21}$ and $z_{22}$ are all distinct and S plays $v_1$.}
    			\label{fig:last-interior}
\end{figure}

By exclusion of \ref{case:parent}, $u$ has no leaf, and by exclusion of \ref{case:single-tail} and \ref{case:grandparent+2}, $u$ has at least two white parents, so by exclusion of \ref{case:triple-tail}, $u$ must have exactly two white parents, which we will call $z_{11}$ and $z_{21}$. By exclusion of \ref{case:zugzwang} and \ref{case:endangered-vertex}, $z_{11}$ has another citrus neighbor $v_1$, with a single further white neighbor $z_{12}$, which is not a leaf. Furthermore, the exclusion of \ref{case:single-tail} implies that $z_{12}$ is a parent. We can also find $v_2$ and $z_{22}$ which lie in analogous relation to $z_{21}$ (see \Cref{fig:last-interior}).

If either $z_{11}=z_{22}$ or $z_{21}=z_{12}$, then D can earn 14 points by playing $u$ as it makes $v_2$ or $v_1$ unplayable. If $z_{11}\neq z_{22}$ and $z_{21}\neq z_{12}$, it follows that $v_1\neq v_2$, but it could be that $z_{12}=z_{22}$. In this case, D should still play $u$, earning 12 points. If S responds by playing $v_1,v_2,z_{12}$ or the leaf of $z_{12}$, that will earn at least 10 points and we are done. If S plays any other move, that will earn 6 points, and D can then play $z_{12}$ for 20 points. Together with the minimum of 6 points from any subsequent move by S, we earn a total of 44 points in four moves.

This leaves the possibility that all of the vertices introduced are distinct (see \Cref{fig:last-interior}). In this situation, D should again play $u$, earning 12 points. If S replies by playing a vertex that is adjacent to both $z_{12}$ and $z_{22}$, then that earns at least 16 points and we are done. Otherwise, we may assume without loss of generality that S plays a vertex that is not adjacent to $z_{22}$, in which case D may earn 16 points by playing $z_{22}$ next. Together with the minimum of 12 points that S has to earn in two moves, that comes to a total of at least 40 points over four moves. 

\subsection*{Phase 10}\label{phase:parents}
Before we go into any detail on the strategy of \nameref{phase:parents}, we present a lemma which we will apply repeatedly. Recall that the game is at time $t$ \emph{after} $t$ moves have been played.

\begin{lemma}\label{lemma:p-q}
Let $P$ be a non-empty set of undominated parents and $Q$ a set of citrus vertices with the property that all undominated neighbors of vertices in $Q$ lie inside $P$. Let further $Y_1\subset Q$ be a set of yellow vertices within $Q$, and let $Y_2$ be a set of yellow vertices, disjoint from $Q$, such that all non-leaf white neighbors of vertices in $Y_2$ lie in $P$. Suppose that the game is at time $t$ and that D is to play next.

If $3\lv P\rv\leq 4\lv Q\rv + \lv Y_1\rv+ \lv Y_2\rv$ or $\lv P\rv$ is odd and $3\lv P\rv\leq 4\lv Q\rv + \lv Y_1\rv+ \lv Y_2\rv+1$, then D can play and we can color in a way such that the game will arrive at a time $t'>t$, such that D is next to play and at least $10(t'-t)$ points have been earned between times $t$ and $t'$.
\end{lemma}

We will apply this lemma several times to sets of explicitly listed vertices, and when listing the elements of $Q, Y_1$ and $Y_2$, we will have to be careful not to include the same vertex twice by a different name. Note that for $P$ this never poses a problem, even when the parity of $\lv P \rv$ changes. We may omit $Y_1$ and $Y_2$ when saying how the lemma is to be applied, meaning that we take them to be the empty set.

\begin{proof}[Proof of \Cref{lemma:p-q}.]
Dominator's strategy is to simply play undominated parents in $P$ until none are left. Suppose this is the case and D is to move next at time $t'>t$. We abbreviate $T=t'-t$, $p=\lv P\rv,q=\lv Q\rv, r_1=\lv Y_1\rv$ and $r_2=\lv Y_2\rv$. Furthermore, we consider the following quantities.
\begin{itemize}
\item $m_D$ is the number of moves made by D between times $t$ and $t'$.
\item $m_S$ is the number of moves made by S in this interval.
\item $\alpha$ is the number of vertices in $Q$ played by S in this interval.
\item $\beta$ is the number of \emph{citrus} vertices played by S in this interval that are \emph{not} in $Q$.
\end{itemize}
Note that we must have $m_S \leq m_D$, $T=m_D+m_S$ and $\alpha \leq p-m_D$. 

We now color vertices and establish a lower bound on the number of points earned since $t$. Whenever D plays a white vertex in $P$, we can color that vertex and its leaf red, which earns $12m_D$ points. Secondly, after move $t'$, no vertex in $Q$ has an undominated neighbor, whence we may color all of them red, which earns $4q+r_1$ points. Next, we know that all vertices in $P$ that D has not played have been dominated by S, and we can color them green, earning an additional $2(p-m_D)$ points. We also know that the vertices in $Y_2$ are not adjacent to any undominated parent anymore, so we can color them green to earn $r_2$ more points.

Now we look at the set of citrus vertices that S has played outside of $Q$. We color each of them red, which earns at least $4\beta$ points (it poses no problem that this may include vertices from $P$ or $Y_2$, which we have just colored green). Note that S must also dominate a new vertex with each such move, and we color all of them green unless we have already done so because they lie in $P$. Since S can have dominated at most $p-m_D-\alpha$ vertices in $P$ by playing vertices outside of $Q$, this earns at least $2(\beta-(p-m_D-\alpha))$ points. Finally, whenever S plays a white vertex outside $P$ or an orange vertex next to a white vertex, we can color the involved white vertex red. This earns $6(m_S-\beta-\alpha)$ points.

We obtain the following lower bound on the number $e$ of points earned.
\begin{align}\label{eq:earned-base}
e&\geq 12m_D+4q+r_1+2(p-m_D)+r_2+4\beta+2(\beta-(p-m_D-\alpha))+6(m_S-\beta-\alpha).
\end{align}
Letting $b=1$ if $p$ is odd and $b=0$ if $p$ is even, we can simplify this to yield
\begin{align*}
e&\geq 12m_D+4q+r_1+r_2+6m_S-4\alpha\\
&\geq 12m_D+4q+r_1+r_2+6m_S+4(m_D-p)\\
&\geq 10T-p-b+6m_D-4m_S.
\end{align*}
If $e$ were less than $10T$, then the above inequality would imply that $p>6m_D-4m_S-b$, and therefore $p\geq 6m_D-4m_S+b$. In this case, we can make use of the fact that all terms on the right-hand side of \eqref{eq:earned-base} must be non-negative and obtain 
\begin{align*}
e&\geq  12m_D+4q+r_1+r_2+2(p-m_D)\\
&\geq 10m_D+5p-b\\
&\geq 10m_D+5(6m_D-4m_S+b)-b\\
&\geq 10T,
\end{align*}
which is in contradiction to the assumption that $e<10T$.
\end{proof}

We now turn to the description of \nameref{phase:parents}. As in the previous phase, we loop through a list of cases, and we continue to recolor vertices such that the number of points associated with their color is as low as possible ahead of each run of the loop. One consequence of the recoloring is that by the end condition of \nameref{phase:white-interior}, no blue vertices exist anymore. The phase ends when none of the cases applies. The list still includes \ref{case:endangered-vertex}, \ref{case:double-tail+}  and \ref{case:zugzwang} from \nameref{phase:white-interior}, and we add the following two further cases.

\begin{enumerate}[label=Case \arabic*., ref=Case \arabic*]
\setcounter{enumi}{\value{case_enum}}
\item\label{case:p-q} Wherever \Cref{lemma:p-q} is applicable, D plays and we color according to its instructions.
\item\label{case:four-parents} Suppose there exists a citrus vertex $u$ with four white parents $x_1,x_2,x_3,x_4$.

Let $i\in \{1,2\}$ and observe that by exclusion of \ref{case:endangered-vertex} and \ref{case:zugzwang}, each $x_i$ has a citrus neighbor $v_i$ with just one additional white neighbor $y_i$ that is not a leaf.\footnote{The same is true for $i\in \{3,4\}$ of course, but we do not make use of that here.} By exclusion of \ref{case:single-tail}, each $y_i$ must be a parent (see \Cref{fig:four-parents}). D can now play $u$. If S replies by playing a (citrus) vertex $z$ that dominates both $y_1$ and $y_2$, we color $z$ and $u$ red and $x_1,\ldots,x_4,y_1$ and $y_2$ green to earn at least 20 points.

If that is not the case, then we may assume by symmetry that S plays a vertex $z$ that is not adjacent to $y_2$, and D can play $y_2$ themselves. We now color $u,z,v_2,y_2$ and the leaf of $y_2$ red and $x_1,x_2,x_3,x_4$ and all white neighbors of $z$ green, thus earning at least 34 points. Together with 6 more points from the subsequent move by S, at least 40 points are earned over four moves.

\begin{figure}[htbp]\centering
    			\includegraphics{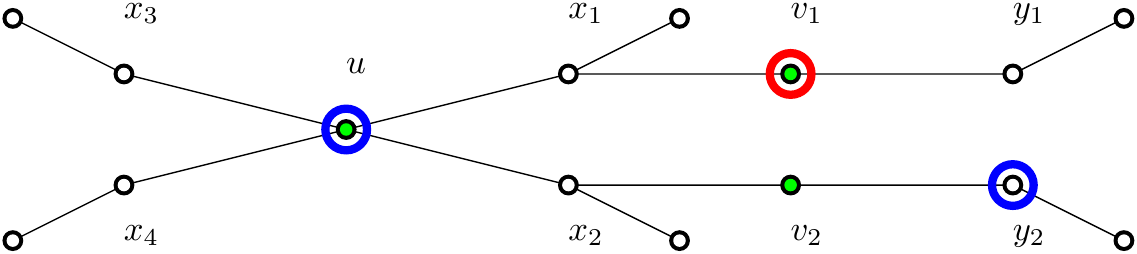}
    			\caption{The situation in which S plays $z=v_1$.}
    			\label{fig:four-parents}
\end{figure}
\end{enumerate}

If none of the above cases apply, there can be no white parent left. We establish this in a series of observations.

\begin{observation}\label{obs:no-two-parents-and-leaf}
If a parent has a white leaf, then it cannot be adjacent to two white parents.
\end{observation}

\begin{proof}
By exclusion of \ref{case:double-tail+}.
\end{proof}

\begin{observation}\label{obs:B-no-green}
If a vertex $u\in B$ is adjacent to a white parent $x_1$ that is not quasi-dependent, then $u$ cannot be green.
\end{observation}
\begin{proof}
Suppose towards a contradiction that a green vertex $u\in B$ were adjacent to a white parent $x_1$ that is not quasi-dependent. Since $u\in B$, $x_1$ must be in $B$ by \Cref{obs:undominated-parent-in-A}. By \Cref{obs:yellow-cousin}, $u$ cannot have a white leaf and by \Cref{obs:yellow-bad-uncle}, $u$ cannot be adjacent to three white parents. On the other hand, by exclusion of \ref{case:endangered-vertex}, we know that $u$ must be adjacent one other white parent $x_2\in B$. Again by \Cref{obs:yellow-bad-uncle}, $u$ cannot be a mixed quasi-anchor, so $x_2$ cannot be quasi-dependent. Because $x_1$ and $x_2$ are both in $B$ and not quasi-dependent, they must each have an additional neighbor in $B$, which we call $v_1$ and $v_2$, respectively. Note that it is possible that $v_1=v_2$.

Regardless of whether this is the case, $v_1$ and $v_2$ must both be yellow by \Cref{obs:core-parents}. If $v\in \{v_1,v_2\}$ has no white leaf, then $v$ must be adjacent to exactly three white parents, two of which are not in $\{x_1,x_2\}$. Indeed, by exclusion of \ref{case:four-parents}, $v$ cannot be adjacent to more than three white parents, and if $v$ were adjacent to only one white parent $y$ that is distinct from $x_1$ and $x_2$ (or to no such parent), we could apply \Cref{lemma:p-q} to $P=\{x_1,x_2,y\}$ (or to $P=\{x_1,x_2\}$, respectively) and $Q=\{u,v\}$. On the other hand, if both $v_1$ and $v_2$ had a white leaf, then they would have to be distinct by \Cref{obs:no-two-parents-and-leaf}, and we could apply \Cref{lemma:p-q} to $P=\{x_1,x_2\}$, $Q=\{u\}$ and $Y_2=\{v_1,v_2\}$. 

We may therefore assume without loss of generality that $v_1$ has no white leaf and is thus adjacent to two white parents $y_{11}$ and $y_{12}$ that are distinct from $x_2$. By exclusion of \ref{case:endangered-vertex} and \ref{case:zugzwang}, both parents must have a neighbor without a white leaf and being adjacent to exactly one additional white parent. We denote these neighbors of $y_{11}$ and $y_{12}$ by $w_{11}$ and $w_{12}$, respectively. Note that $w_{11}$ and $w_{12}$ can be the same vertex but must be different from $u,v_1$ and $v_2$. If $w_{11}=w_{12}$, then we can apply \Cref{lemma:p-q} to $P=\{x_1,x_2,y_{11},y_{12}\}$ and $Q=\{u,v_1,w_{11}\}$. Otherwise, we denote the respective additional white parents of $w_{11}$ and $w_{12}$ by $z_{11}$ and $z_{12}$, respectively. 

If in this situation $v_2$ has a white leaf (and is therefore not adjacent to any white parent other than $x_2$), then we can apply \Cref{lemma:p-q} to $P=\{x_1,x_2,y_{11},y_{12},z_{11},z_{12}\}$, $Q=\{u,v_1,w_{11},w_{12}\}$, $Y_1=\{v_1\}$ and $Y_2=\{v_2\}$. If $v_2$ has no white leaf, then without loss of generality, we can find vertices $y_{21},y_{22},w_{21},w_{22},z_{21}$ and $z_{22}$ in a configuration that is identical to the configuration formed by $y_{11},y_{12},w_{11},w_{12},z_{11}$, except that $y_{21}$ and $y_{22}$ are adjacent to $v_2$ instead of $v_1$. Because $x_2\notin \{x_1,y_{11},y_{12}\}$, we must have $v_1\neq v_2$, and by considering the number of adjacent white parents of citrus vertices, it is easy to see that $\{v_1,v_2\}$ and $\{w_{11},w_{12},w_{21},w_{22}\}$ are disjoint.

\begin{figure}[htbp]\centering
    			\includegraphics{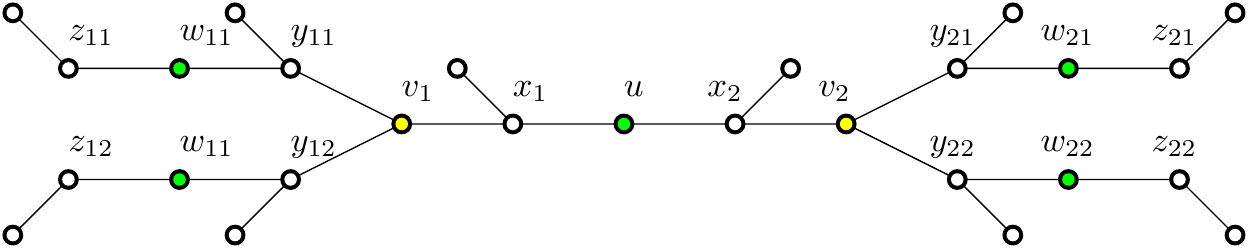}
    			\caption{The situation where all introduced vertices are distinct.}
    			\label{fig:green-in-B}
\end{figure}

If $\{w_{11},w_{12}\}$ and $\{w_{21},w_{22}\}$ are also disjoint (see \Cref{fig:green-in-B}), we can apply \Cref{lemma:p-q} with 
\begin{align*}
P&=\{x_1,x_2,y_{11},y_{12},z_{11},z_{12},y_{21},y_{22},z_{21},z_{22}\},\\
Q&=\{u,v_1,v_2,w_{11},w_{12},w_{21},w_{22}\} \text{ and }\\
Y_1&=\{v_1,v_2\}.
\end{align*}
If $\{w_{11},w_{12}\}$ and $\{w_{21},w_{22}\}$ are not disjoint, it is easy to check that \Cref{lemma:p-q} is still applicable for the same sets $P$ and $Q$. For example, if $w_{11}=w_{21}$, then we must have $\{y_{11}, z_{11}\}=\{y_{21}, z_{21}\}$ also.
\end{proof}

\begin{observation}\label{obs:cousins-resolved}
If two adjacent parents each have an undominated leaf, both parents must themselves be dominated.
\end{observation}
\begin{proof}
Suppose there exist parents $u$ and $x$ of white leaves for which the claim is not true because $x$ is not dominated. By \Cref{obs:no-cousins-inside-A}, both $u$ and $x$ must lie in $B$, and by \Cref{obs:quasi-anchor-no-leaf}, neither $u$ nor $x$ can be quasi-dependent. Therefore, $x$ must have another neighbor $v\in B$. By \Cref{obs:B-no-green}, both $u$ and $v$ must be yellow. Also, $u$ has no white parent other than $x$ by \Cref{obs:no-two-parents-and-leaf}.

If $v$ has a white leaf, then $v$ has no other white parents either, and we can apply \Cref{lemma:p-q} to $P=\{x\}$, $Q=\emptyset$ and $Y_2=\{u,v\}$. If $v$ does not have a white leaf, we know that by exclusion of \ref{case:endangered-vertex} and \ref{case:four-parents}, $v$ must have one or two white neighbors other than $x$. If $v$ is adjacent to exactly one additional white parent $y_1$, we can apply \Cref{lemma:p-q} to $P=\{x,y_1\}$, $Q=\{v\}$, $Y_1=\{v\}$ and $Y_2=\{u\}$.

This leaves the possibility that $v$ has exactly two additional neighbors $y_1$ and $y_2$. By exclusion of \ref{case:zugzwang} and \ref{case:endangered-vertex}, $y_1$ and $y_2$ have respective neighbors $w_1$ and $w_2$, each having no white leaf and being adjacent to exactly two white parents. Note that both of them must be different from $u$ and $v$. If $w_1=w_2$, we can apply \Cref{lemma:p-q} with $P=\{x,y_1,y_2\}$ and $Q=\{v,w_1\}$.  Otherwise, we denote the additional white parents adjacent to $w_1$ and $w_2$ by $z_1$ and $z_2$, respectively, and apply \Cref{lemma:p-q} to $P=\{x,y_1,y_2,z_1,z_2\}$, $Q=\{v,w_1,w_2\}$, $Y_1=\{v\}$ and $Y_2=\{u\}$.
\end{proof}

\begin{observation}\label{obs:double-yellow}
No white parent can have two yellow neighbors.
\end{observation}
\begin{proof}
The proof of \Cref{obs:double-yellow} is similar to but simpler than those of \Cref{obs:B-no-green} and \Cref{obs:cousins-resolved} and we leave it to the reader.
\end{proof}

\begin{observation}\label{obs:semi-tridents-resolved}
No semi-trident is adjacent to any white parent other than its anchor-bridges.
\end{observation}
\begin{proof}
Let $v$ be a semi-trident that is adjacent to a white parent $x$ which is not an anchor-bridge. By \Cref{obs:semi-trident-residual}, $x$ is not residual, so $x$ must have a neighbor $u\in B$. By \Cref{obs:cousins-resolved} and exclusion of \ref{case:endangered-vertex}, $u$ is adjacent to another white parent $y$, which cannot be in $A$ by \Cref{obs:undominated-parent-in-A}.

Now we invoke \Cref{obs:semi-trident-transitional} to see that at most one of $x$ and $y$ can be quasi-dependent. We may thus assume that $y$ is not quasi-dependent, and therefore has another neighbor $w\in B$. By \Cref{obs:B-no-green}, both $u$ and $w$ have to be yellow, which contradicts \Cref{obs:double-yellow}.
\end{proof}

To make our final two observations, we need to introduce some more terminology. We say that an anchor is \emph{true} if it still has a white dependent parent and that a white anchor-bridge is \emph{true} if it is adjacent to a true anchor. The \emph{split} of a white parent is the number of its non-leaf neighbors, i.e., its degree minus one. If a vertex is adjacent to a white parent but is not a leaf, then we call it an \emph{uncle}. Note that uncles cannot themselves have white leaves by \Cref{obs:cousins-resolved}. The \emph{valency} of an uncle is the number of its white neighbors, which must all be parents.

\begin{observation}\label{obs:three-on-two-cororllary}
White parents and their neighbors have the following properties.
\begin{enumerate}[label=(\roman*)]
\item The valency of an uncle must be 2 or 3.
\item Every white parent is adjacent to exactly one uncle of valency 2.
\item A true anchor has valency 2.
\item No parent can be adjacent to two or more true anchors.
\item A true anchor-bridge has split 2.
\item No uncle is adjacent to two true anchor-bridges.
\end{enumerate}
\end{observation}
\begin{proof}
We prove these properties one-by-one.
\begin{enumerate}[label=(\roman*)]
\item By exclusion of \ref{case:endangered-vertex}, the valency cannot be 1, and by exclusion of \ref{case:four-parents}, it cannot be more than 3. 
\item It follows by exclusion of \ref{case:zugzwang} that every white parent must have at least one neighbor of valency 2. Suppose now that a white parent $x$ has two neighbors $u$ and $v$ of valency 2. Letting $y$ and $z$ be the second white parents adjacent to $u$ and $v$, respectively, we can apply \Cref{lemma:p-q} to $P=\{x,y,z\}$ and $Q=\{u,v\}$.
\item A true anchor has at least one white parent by definition. If it were adjacent to more than two white parents, then \ref{case:zugzwang} would apply to its white dependent parent.
\item Follows from (ii) and (iii).
\item A true anchor-bridge must have split 2 because any anchor-bridge has degree 3.
\item By \Cref{obs:anchor-bridge-uncles}, no uncle can have more than two true anchor-bridges. A vertex with exactly two true anchor-bridges would be a semi-trident and would have valency 2 by \Cref{obs:semi-tridents-resolved}. But then both of the true anchor-bridges would violate (ii) since the true anchors they are adjacent to also have valency 2.\qedhere
\end{enumerate}
\end{proof}

\begin{observation}\label{obs:parents-resolved}
There is no white parent left.
\end{observation}
\begin{proof}
We want to apply \Cref{lemma:p-q} with $P$ being the set of all white parents and $Q$ being the set of all uncles, but we have to find a suitable upper bound on $\lv P\rv$ in terms of $\lv Q\rv$ first. For this purpose, consider the following quantities.

\begin{enumerate}
\item[$\alpha$:] The number of white dependent parents.
\item[$\beta$:] The number of true anchor-bridges.
\item[$\gamma$:] The number of white parents with a split of 2 that are not true anchor-bridges, i.e., are not adjacent to a true anchor.
\item[$\delta$:] The number of white parents with a split of 3 or higher that are adjacent to a true anchor.
\item[$\varepsilon$:] The number of white parents with a split of 3 or higher that are not adjacent to a true anchor.
\item[$\phi$:] The number of true anchors.
\item[$\chi$:] The number of uncles that are not true anchors but have valency 2.
\item[$\psi$:] The number of uncles of valency 3 that are adjacent to a true anchor-bridge.
\item[$\omega$:] The number of uncles of valency 3 that are not adjacent to a true anchor-bridge.
\end{enumerate}

Note that any white parent or non-leaf neighbor of a white parent contributes to exactly one of the above quantities by \Cref{obs:three-on-two-cororllary}.(i). Therefore, saying that a vertex is of \emph{type} $\alpha, \beta$ etc. is well-defined. Note that $\lv P\rv = \alpha + \beta + \gamma + \delta + \varepsilon$ and $\lv Q\rv = \phi + \chi+\psi+\omega$.

We can state a number of linear constraints between these quantities. To begin with, every true anchor is by \Cref{obs:no-two-dependent-parents} and \Cref{obs:three-on-two-cororllary}.(iii) adjacent to exactly one  dependent and one non-dependent parent white parent. The non-dependent parent must be either of type $\beta$ or of type $\delta$ and by \Cref{obs:three-on-two-cororllary}.(iv), no parent is adjacent to more than one anchor. Therefore, we have 
\begin{align}\label{eq:final1}
\alpha=\phi= \beta + \delta.
\end{align}
On the other hand, every true anchor-bridge must be adjacent to a vertex of valency 3 by \Cref{obs:three-on-two-cororllary}.(ii), whence by \Cref{obs:three-on-two-cororllary}.(vi),
\begin{align}\label{eq:final2}
\beta\leq \psi.
\end{align}
Next, every white parent that is not adjacent to an anchor must be adjacent to some other uncle with valency 2, so we have 
\begin{align}\label{eq:final3}
\gamma+\varepsilon\leq 2\chi.
\end{align} 
Lastly, every parent of type $\delta$ or $\varepsilon$ is adjacent to two uncles of valency 3, and every parent of type $\beta$ or $\gamma$ is adjacent to one uncle of valency 3. Therefore, 
\begin{align}\label{eq:final4}
\beta+\gamma+2\delta+2\varepsilon\leq 3 \psi+3\omega.
\end{align}
Now we can bound $3\lv P\rv$ by $4\lv Q\rv$ as follows.
\begin{align*}
3(\alpha+\beta+\gamma+\delta+\varepsilon)&\overset{\eqref{eq:final1}}{=}  2\beta+3\gamma+2\delta+3\varepsilon+4\phi\\
&\overset{\eqref{eq:final2}}{\leq} \beta+3\gamma+2\delta+3\varepsilon+4\phi+\psi\\
&\overset{\eqref{eq:final3}}{\leq}\beta + \gamma + 2\delta + \varepsilon + 4\phi + 4\chi + \psi\\
&\overset{\eqref{eq:final4}} \leq 4\phi + 4\chi + 4\psi + 3\omega.
\end{align*}
\end{proof}

\subsection*{Phase 11}\label{phase:isolates}
This phase lasts until there is no undominated vertex left, i.e., until the game ends. It is sufficient to show that any move earns at least 10 points. This is the case as every white vertex left is either a leaf with a green parent or has two orange neighbors.

This concludes the proof of \Cref{prop:main}.

\section{Concluding remarks}

\subsection{The complexity of finding the next move}

The proof of \Cref{prop:main} gives algorithmic instructions for every move of Dominator, and the algorithm could be implemented as computer software if desired. What is more, when playing on an isolate-free graph $G$ with $m$ edges, a slightly modified algorithm could compute each of Dominator's moves in only $O(m)$ operations. We sketch briefly how this may be achieved.

Given a list of all edges, we may obtain the neighborhoods of all vertices in $O(m)$ steps. By running thrice through all vertices and their neighborhoods, we can find all leaves, dependent parents and anchors. Using this data, we prune $G$ according to the proof of \Cref{lemma:structure} until $G$ has the properties listed in the lemma. Each pruning step removes at least one edge from $G$ while adding a uniformly bounded number of operations as a prefix to the calculation of any move, thus not jeopardizing the bound $O(m)$.

The end conditions of all phases except \nameref{phase:big-white-components} and \nameref{phase:parents} are formulated as conditions on local configurations that can be checked by running a fixed number of times through all vertices and their neighborhoods, marking those with undominated leaves, dependent parents etc. and then considering the white neighbors of each vertex. If the end condition is not met, the local configuration witnessing this will always allow us to find an active move or a sequence of moves that D may play. Likewise, if S makes a move that requires a specific reaction, the correct reactive move can also be found quickly. 

For \nameref{phase:big-white-components}, we only need to reformulate the end condition as stipulating that no white vertex has two white neighbors. Observe that a white cycle or longest white path can be found in $O(m)$ steps if no white vertex has three white neighbors.

This leaves \nameref{phase:parents}. Its end condition looks as if the algorithm would need to check whether \ref{case:p-q} applies to any combination of adequate sets, but in reality it is sufficient to find a configuration of vertices that does not satisfy one the subsequent observations. The proof of the first such observation will then reference a case among \ref{case:endangered-vertex},  \ref{case:double-tail+}, \ref{case:zugzwang}, \ref{case:p-q}  and \ref{case:four-parents} with instructions that can be calculated in $O(m)$.

\subsection{Extremal graphs for \Cref{thm:main}}
\Cref{thm:main} is tight for both the Dominator-start and the Staller-start version of the game and there is a large family of examples that show this. Indeed, for any graph $G$, we can create an extremal graph $\hat{G}$ for the Dominator-start domination game by appending two dependent parents to each vertex. During the game on $\hat{G}$, whenever D plays one of the appended parents or their leaf, S plays the vertex in $G$ they are appended to. If D plays a vertex from $G$ and there is still another vertex in $G$ left, then S plays that. This way, every vertex in $G$ and one vertex of every appended parent-leaf pair will end up being played, which makes three out of five vertices in $\hat{G}$.

The family of extremal examples for both versions of the game can be further extended by taking any existing extremal example $G$, selecting any dependent parent $v$ in $G$, and generating a new graph $G_v$ by removing the leaf of $v$ and appending three dependent parents to $v$. We denote by $\mathcal{C}$ the family of all graphs generated iteratively this way, starting from some $\hat{G}$.\footnote{In the proof of \Cref{thm:main}, this procedure will always be reversed by \Cref{lemma:structure}.(iii).} This construction is a straightforward generalization of a similar construction for trees found by Henning and Löwenstein in \cite{henning_extremal_trees}.

A graph of the form $\hat{G}$ can be modified further into an extremal example $\tilde{G}$ for the Staller-start version by appending a single leaf to an arbitrarily chosen vertex in $V(G)\subset V(\hat{G})$. When playing the game in $\tilde{G}$, S plays this leaf on the first move and then plays as in $\hat{G}$ for the rest of the game.


It is worth asking to what extent these constructions exhaust the extremal graphs for \Cref{thm:main}. We point out that $C_5$ is also extremal for the Dominator-start version (and we may append a leaf to it to obtain an extremal graph for the Staller-start version), but this could well turn out to be a minor caveat. Analyzing the original version of the game without gifted moves, one cannot even take multiple disconnected copies of $C_5$ without moving far away from the the bound of $3n/5$. This leads to the following question.

\begin{question}\label{question:copy-robust}
Does there exist a non-empty graph $G\notin \mathcal{C}$ on $n$ vertices so that for all $k\in \N$ the graph $kG$ consisting of $k$ disconnected copies of $G$ satisfies $\gamma_g(kG)=3kn/5$?
\end{question}

Even more interesting than an extremal graph outside of $\mathcal{C}$ that is robust with respect to disjoint copies would be a large connected example.

\begin{question}\label{question:connected}
Does there exist, for every $N\in \N$, a connected graph $G\notin \mathcal{C}$ on $n>N$ vertices with $\gamma_g(G)=3n/5$?
\end{question}

Of course, analogous questions could be asked about the Staller-start version of the game, but it appears that they would be largely equivalent to \Cref{question:copy-robust} and \Cref{question:connected}. Indeed, whenever we find a graph $G$ such that $\gamma_g'(G)=(3n+2)/5$, we may consider Staller's first move $x$ and note that the graph $H$ consisting of $G$ with $x$ and all its leaves removed must satisfy $\gamma_g(H)=3(n-1)/5$. In the other direction, if a reasonably large graph $G$ satisfies $\gamma_g(G)=3n/5$, it seems likely that we can find a suitable vertex to which we can append a leaf which Staller can play as their first move.

\subsection{Unifying different variants of the domination game}
Research on the domination game has led to the introduction of several new variants of the game. Most prominently, the \emph{total domination game}, first studied in \cite{henning_total_initial}, is a variant in which every vertex must have a neighbor in the set of selected vertices for the game to end, regardless of whether the vertex has been selected itself. Analogously to $\gamma_g$, one defines the \emph{game total domination number} $\gamma_{tg}$. 

For the total domination game, it is natural to consider only graphs that do not have isolated vertices or edges, and in \cite{henning_total_4/5_and_conjecture}, Henning,  Klav{\v{z}}ar, and Rall showed that for such graphs $G$, the bound $\gamma_{tg}(G)\leq 4n/5$ holds. Furthermore, they conjectured that $\gamma_{tg}(G)\leq 3n/4$. Following further progress by Bujt{\'a}s \cite{bujtas_total_11/14}, who showed that $\gamma_{tg}(G)\leq 11n/14$, Portier and the author have recently confirmed the $3/4$-conjecture \cite{portiertotaldomination}.
\begin{theorem}\label{thm:total}
For a graph $G$ on $n$ vertices without isolated vertices or edges, $\gamma_{tg}(G)\leq 3n/4$.
\end{theorem}
Another variant, the \emph{hypergraph transversal game}, was introduced by Bujt{\'a}s, Henning, and Tuza in \cite{bujtas_transversal}. This game is played on a hypergraph $H=(V_H,E_H)$ by two players, who select vertices alternately until the set of selected vertices $A$ becomes \emph{transversal}, i.e., until it contains a vertex of every edge in $E_H$. Of course, one can define the \emph{game transversal number} $\tau_g$ analogously to $\gamma_g$ and $\gamma_{tg}$. Bujt{\'a}s et al. proved the following theorem about $\tau_g$, which we state here in simplified form. 
\begin{theorem}\label{thm:transversal}
Let $H$ be a hypergraph with $n_H$ vertices and $m_H$ (hyper-)edges. If all edges in $H$ have size at least two and $H$ does not belong to a small exceptional class of hypergraphs, then
\begin{align*}
\tau_g(H)\leq \frac{4}{11}(n_H+m_H).
\end{align*}
\end{theorem}
This result is certainly interesting in its own right, but it can also be used to derive bounds on $\gamma_g$ and $\gamma_{tg}$. For instance, in \cite{bujtas_transversal} it is shown that for any graph $G$ with minimum degree $\delta(G)\geq 2$, $\gamma_{tg}(G)\leq 8n/11$. This is done by applying \Cref{thm:transversal} to the hypergraph $H_G$ that has $V(G)$ as its vertex set and the neighborhoods of vertices of $G$ as its edges, noting that $\tau(H_G)=\gamma_{tg}(G)$. Similarly, bounds on $\gamma_g(G)$ can be derived by considering the hypergraph $\overline{H_G}$ whose edges are the closed neighborhoods of vertices in $G$.

The specific bound $\gamma_{tg}(G)\leq 8n/11$ for graphs $G$ with minimum degree at least 2 has since been superseded \cite{portiertotaldomination}, but the general approach of viewing variants of the domination game as special cases of the hypergraph transversal game still looks promising. Indeed, it seems likely that by restricting the class of hypergraphs under consideration, one can improve on \Cref{thm:transversal} substantially, without sacrificing the ability to deduce results about the domination game or variants thereof. For instance, we know that for any graph $G$ without isolated vertices, every vertex must lie in at least one edge of $H_G$, $H_G$ must satisfy $m_H\leq n_H$, and if two edges of $H_G$ share a vertex, then there is a third edge intersecting both of the other edges. But \Cref{thm:transversal} assumes none of these properties. 

That being said, it might also be the case that bounds derived from the game transversal number will ultimately always be weaker than those obtained from arguments about the specific game. The following question is intended to put the power of the approach to the test.

\begin{question}
Can \Cref{thm:main} or \Cref{thm:total} be derived from a common theorem\footnote{Of course, one can trivially reformulate \Cref{thm:main} and \Cref{thm:total} as a two-parted statement about the classes of hypergraphs of the form $H_G$ and $\overline{H_G}$, but that is not in the spirit of the question.} about the hypergraph transversal game?
\end{question}

We believe that even a partial answer to this question could move the area forward by revealing deeper insights about the relationships between different variants of the domination game, which as of now seem somewhat isolated. An affirmative answer in particular should lead to progress on other open questions in the area.

\printbibliography

\end{document}